\documentclass[fleqn,12pt,a4paper]{wlscirep}
\usepackage[utf8]{inputenc}
\usepackage[T1]{fontenc}
\usepackage{url,hyperref,lineno,microtype,subcaption}
\usepackage[onehalfspacing]{setspace}
\usepackage{longtable}
\usepackage{float}
\usepackage{amsmath}
\usepackage{color}
\usepackage{amssymb,amsthm}
\usepackage[ruled, vlined]{algorithm2e}
\usepackage{algorithmic}

\def\d{\delta}

\newtheorem{theorem}{Theorem}[section]

\newtheorem{definition}[theorem]{Definition}

\title{Spatiotemporal slope stability analytics for failure estimation (SSSAFE): linking radar data to the fundamental dynamics of granular failure}
\author[1,*]{Antoinette Tordesillas}
\author[1]{Sanath Kahagalage}
\author[2]{Lachlan Campbell}
\author[2]{ \\ Pat Bellett}
\author[3]{ Emanuele Intrieri }
\author[4]{ Robin Batterham }
\affil[1]{School of Mathematics $\&$ Statistics, University of Melbourne, Australia}
\affil[2]{GroundProbe, Orica, Australia}
\affil[3]{Department of Earth Sciences, University of Florence, Italy}
\affil[4]{Melbourne School of Engineering, University of Melbourne, Australia}

\affil[*]{atordesi@unimelb.edu.au}


\begin{abstract} 
Impending catastrophic failure of granular earth slopes manifests distinct kinematic patterns in space and time.  While risk assessments of slope failure hazards have routinely relied on the monitoring of ground motion, such precursory failure patterns remain poorly understood. A key challenge is the multiplicity of spatiotemporal scales and dynamical regimes.  In particular, there exist a precursory failure regime where two mesoscale mechanisms coevolve, namely, the preferred transmission paths for force and damage.  Despite extensive studies, a formulation which can address their coevolution not just in laboratory tests but also in large, uncontrolled field environments has proved elusive. Here we address this problem by developing a slope stability analytics framework which uses network flow theory and mesoscience to model this coevolution and predict emergent kinematic clusters solely from surface ground motion data.  We test this framework on four data sets: one at the laboratory scale using individual grain displacement data; three at the field scale using line-of-sight displacement of a slope surface, from ground-based radar in two mines and from space-borne radar for the 2017 Xinmo landslide. The dynamics of the kinematic clusters deliver an early prediction of the geometry, location and time of failure. 
\end{abstract}

\begin{document}

\flushbottom
\maketitle

\thispagestyle{empty}


\section*{Introduction}
\label{Introduction}

Natural and engineered slopes, composed of granular materials like rocks, concrete and soil, can maintain their structural integrity even as damage spreads.  But there is a tipping point, beyond which damage can propagate to cause catastrophic failure with little to no apparent warning signs at the macroscale~\cite{carla2019,Handwerger2019rainfall,Clarkson2020damfailure}.  The landslides in Xinmo (China, 2017) and the dam collapse in Brumadinho (Brazil, 2019) are recent reminders of the devastating impact of slope failure on human lives and livelihoods, infrastructure, and the environment~\cite{carla2019,Handwerger2019rainfall,Clarkson2020damfailure,WLF5,SILVAROTTA2019Brumadinho}. Here the term {\it failure} is used from an operative point of view and is the moment when the slope totally or partially collapses, displaying a paroxysmal acceleration and a disintegration of the mobilized material.  A critical frontline defense against these hazards is large-scale monitoring and analysis of slope movement using remote sensing technologies~\cite{MCQUILLAN2020151,carla2019,Clarkson2020damfailure,Dicketal2015,harries2006,wasowski2014,intrieri2019,Intrieri2017,Dai2020}.  
Some of these measurements have now reached spatial and temporal resolutions (e.g., Slope Stability Radar~\cite{Dicketal2015}) which enable direct connections to be made to the fundamental deformation and failure of granular materials~\cite{MCQUILLAN2020151,Clarkson2020damfailure,carla2019,Intrieri2017}.
Nevertheless, there are significant challenges to overcome before the full potential of these data assets can be harnessed for geotechnical risk assessment and
hazard management~\cite{Intrieri2017,Dai2020}.  One of, if not, the biggest challenge lies in the analysis and interpretation of monitoring data with respect to the underlying micromechanics and dynamics of deformation in the precursory failure regime (PFR)~\cite{intrieri2019,MCQUILLAN2020151,Dicketal2015,wasowski2014,harries2006}.  
Here we address this challenge by formulating a holistic framework for spatiotemporal slope stability analytics for failure estimation (SSSAFE).  SSSAFE is physics-based and bears explicit connections to the micromechanics and dynamics of ductile to brittle failure in granular solids (e.g.,~\cite{ATSKCRMNJT2,ATSKCRMNJT,tordesillas2015,tordesillasEtAl2016} and references therein).

A hallmark of SSSAFE is its detailed characterization of the spatiotemporal coevolution of the preferred pathways for force and damage in PFR using kinematic data.  As highlighted in various reviews~\cite{MCQUILLAN2020151,Clarkson2020damfailure,Dai2020,intrieri2019,Intrieri2017,wasowski2014}, scant attention has been paid to the spatiotemporal dynamics of landslide deformation, with existing approaches in landslide forecasting and early warning systems (EWS) 
falling into one of two categories: (a) spatial analysis of an unstable slope to estimate the location and geometry of a landslide~\cite{hoek1981}, or (b) temporal analysis of ground deformation of single measurement points exhibiting tertiary creep, to deliver a short-term forecast of the time of failure~\cite{Carl2016GuidelinesOT,intrieri2019,carla2019}.  The former partially relies on expert judgment (e.g. the choice of the failure criterion and the method of analysis~\cite{MCQUILLAN2020151}) and on in situ data (depth of the lithologies and of the water table, resistance parameters of the rock or soil) that always bear a certain level of uncertainty and representativeness bias~\cite{ChristianBaecher2011}.

In temporal analysis, the inverse velocity (INV) theory originally proposed by Fukuzono~\cite{fukuzono1985} is the most widely applied method for prediction of the time of collapse in the terminal stages of PFR.  This approach has no spatial aspect and depends on assumptions which motivate areas for improvement in forecasting~\cite{Dicketal2015,carla2019}.   Being based on the inverse value of a derivative parameter, this method is heavily affected by noise, especially when the velocity is not particularly high. While this can be addressed by smoothing the data using a moving average~\cite{Carl2016GuidelinesOT,intrieri2019}, this comes at the cost of diminished sensitivity to important changes in the acceleration trends due to short-terms events, including: surface boundary conditions (e.g., civil engineering and mining works~\cite{MCQUILLAN2020151}), variations in the trigger factors of slope instability (e.g., rainfall~\cite{Handwerger2019rainfall}, seismic~\cite{Qiu2016seismic}, mining blasts~\cite{Dicketal2015,harries2006}), and the inherently complex mechanical interactions between different parts of the slope.  Concurrent sites of instability may also interact and induce stress redistributions that lead a landslide to ``self-stabilize''~\cite{hungr2014,wang2019,TordesillasZhouBatterham2018}. Efforts~\cite{carla2019,Dicketal2015} to improve the INV approach give prima facie evidence to suggest that more accurate forecasts can be achieved when the spatial characteristics of slope displacements are incorporated in the temporal analysis of monitoring data.

Accordingly, recent work focused on the spatiotemporal evolution of landslide kinematics in PFR in two case studies using: (a) ground-based radar data of a rockfall in an open pit mine (Mine 1) where two sites of instability emerged, leading one to self-stabilize before the larger one collapsed; and (b) satellite-based Sentinel 1 radar data (Xinmo) of the catastrophic collapse in Xinmo, which led to 83 fatalities~\cite{TordesillasZhouBatterham2018,KSAT,SDAT,ZhouAOAD,WangSS,WLF5}.  Guided by lessons learned from the physics and dynamics of granular failure, these delivered a reliable early prediction of the location and geometry of the failure region~\cite{TordesillasZhouBatterham2018,ZhouAOAD,WangSS,KSAT,WLF5}, as well as regime change points in PFR~\cite{KSAT,SDAT,WLF5,WangSS}.   In this study, we build on these efforts to develop a holistic data-driven framework which eliminates the uncertainties associated with a postulated stress-strain model for the slope, yet holds explicit connections to the first principles of fracture and failure mechanics of heterogeneous and disordered granular solids (e.g.,~\cite{ATSKCRMNJT2,ATSKCRMNJT,tordesillas2015,tordesillasEtAl2016} and references therein).  To do this, we adopt a transdisciplinary approach which integrates network flow theory of granular failure~\cite{ATSKCRMNJT2,ATSKCRMNJT,tordesillas2015,tordesillasEtAl2016} and mesoscience~\cite{li2013, li2018mesoscience, li2018}.  Given the novelty of this formulation from several fronts, the next section gives a brief review of the relevant developments which, woven together, form the basis of SSSAFE.

\section*{Precursory dynamics of granular failure across system levels and scales}
\label{sec:connect}

\begin{figure}
\centering
\includegraphics[width=1\textwidth]{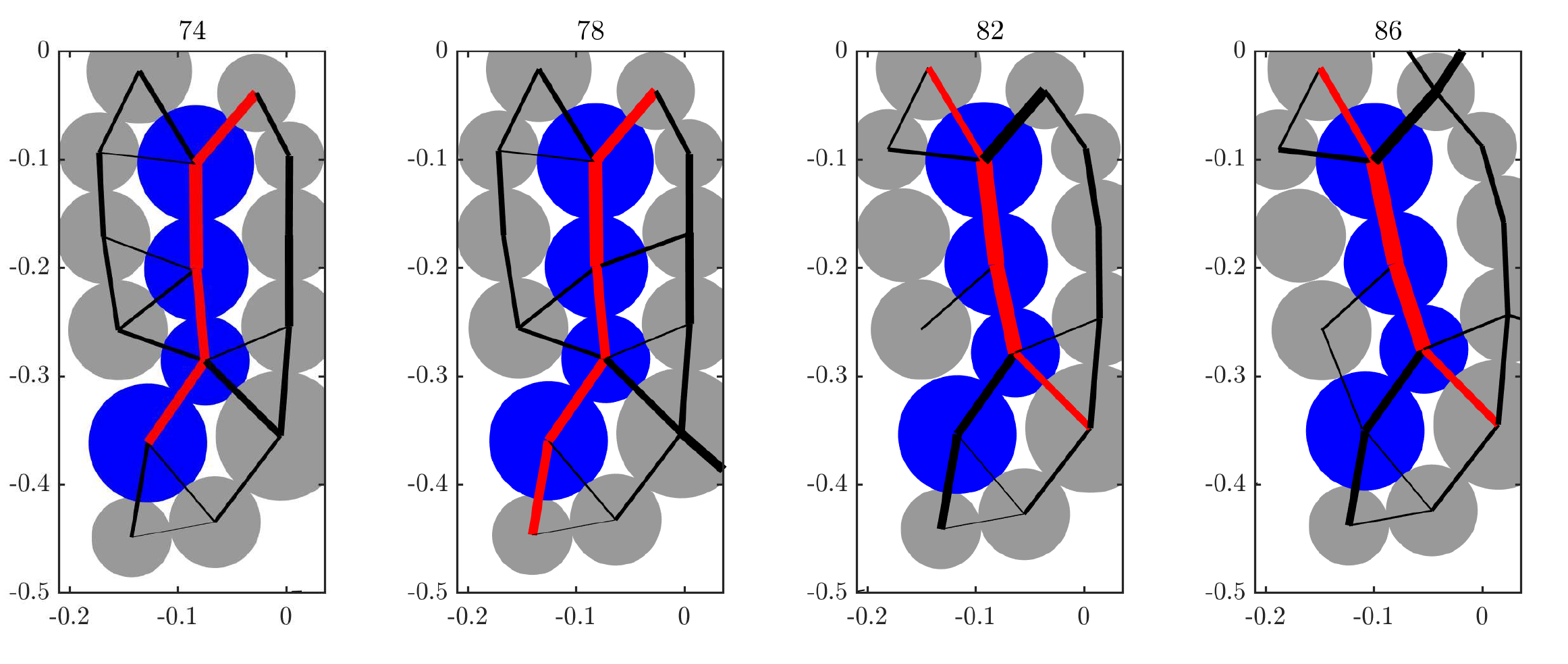}
\caption{(Color online) Redistribution of contact forces around a force chain in the sample Biax across different stages of PFR, prior to the time of failure $t^F_B=104$.  Link thickness is proportional to the contact force magnitude.  Red (black) links correspond to contacts between member particles of the force chain (all other supporting contacts).   Most of the grains in the force chain are colored blue to aid visualization. There is a build up of force across stages 74-78, leading to a new force chain contact at the bottom amid rearrangements of supporting lateral contacts.  Further build up of force in the force chain column results in the buckling of the top and bottom segments of the chain across stages 78-82: in turn, more force is rerouted to the bottom right (top left) in stages 82-86, resulting in a new force chain contact.     }
\label{fig:rewire}
\end{figure}

{\it a) Preferred paths for transmission of damage versus force. \,\,} In complex systems, not all paths for transmission are created equal.  Some are preferred over others.  Experimental studies into the transmission of force and energy in natural and synthetic granular media (e.g., sand, photoelastic disk assemblies) and associated discrete element simulations have shown that the mesoregime of PFR is governed by the coupled evolution of two dominant mechanisms~\cite{burnley2013,ATSKCRMNJT,ATSKCRMNJT2,tordesillas2015,PRE}.  The first comprises the preferred paths for force transmission (mechanism A): a set of system-spanning paths that can transmit the highest force flow along direct and shortest possible routes through the system.  Distinct force chains (Figure~\ref{fig:rewire}) can be readily observed to form along these percolating paths, in alignment with the major principal stress axis~\cite{CundallStrack1979,Potyondy_2004,Majmudar_2005}. The second are the preferred paths for damage (mechanism B), where cracks and/or shearbands emerge. Note that the term {\it damage} is broadly used here to mean the separation of two grains in contact, bonded or unbonded.  

 {\it b) Coevolution of preferred paths: a compromise-in-competition. \,\,}   Arguably the best manifestation of the coupled evolution between force and damage can be observed in deforming photoelastic disk assemblies~\cite{Majmudar_2005,PRE}.  Here one can readily observe forces continually rerouted to alternative pathways as damage spreads (Figure~\ref{fig:rewire}).  This scenario is similar to traffic flows where vehicles are diverted to alternative routes when a road is closed off for repairs or other incidents.  Following this analogy to road networks, grain contact networks similarly give rise to emergent flow bottlenecks.  Prior network flow studies have shown that these sites, which are highly prone to congestion, ultimately become the preferred paths for damage in the failure regime~\cite{ATSKCRMNJT,ATSKCRMNJT2}.  Counter to intuition, however, the bottlenecks do not generally coincide with the location of damage sites in the nascent stages of PFR, which is when ideally predictions should be made to allow enough time to enact mitigative measures.  Instead, a process which can be described as a {\it compromise-in-competition} between the preferred paths for force and damage develops, which effectively shields the bottleneck from damage.  Specifically, force congestion in the bottleneck is relieved by complex stress redistributions that redirect forces to other parts of the sample, where damage can be accommodated with minimal reduction to the system's resistance to failure (or global force transmission capacity).  This may explain why current failure detection methods, which rely on damage sites in the early stages of PFR for spatial clues on where catastrophic failure ultimately forms, sometimes suffer high false positive rates in laboratory~\cite{chakraborty2019early} and field levels~\cite{Intrieri2017}.

\begin{table}
\caption{Mechanisms underlying the strength and failure of granular materials compete in the precursory failure regime (PFR).}
\centering
\begin{tabular}{ccc}
\toprule
 \textbf{Regime}	& \textbf{A -Preferred paths for force }	& \textbf{B -Preferred paths for damage}\\
 \textbf{(emergent structures)}	& \textbf{(force chains)}	& \textbf{(cracks, shear bands)}\\
 \midrule
(A) Stable regime	& dominant 			& suppressed\\
(A-B) Mesoregime PFR  	& compromise-in-competition			& compromise-in-competition\\
(B) Failure regime	& suppressed		& dominant\\
 \bottomrule
 \end{tabular}
 \label{Tab:compete}
 \end{table}

\begin{figure}
\centering
\includegraphics[width=0.9\textwidth]{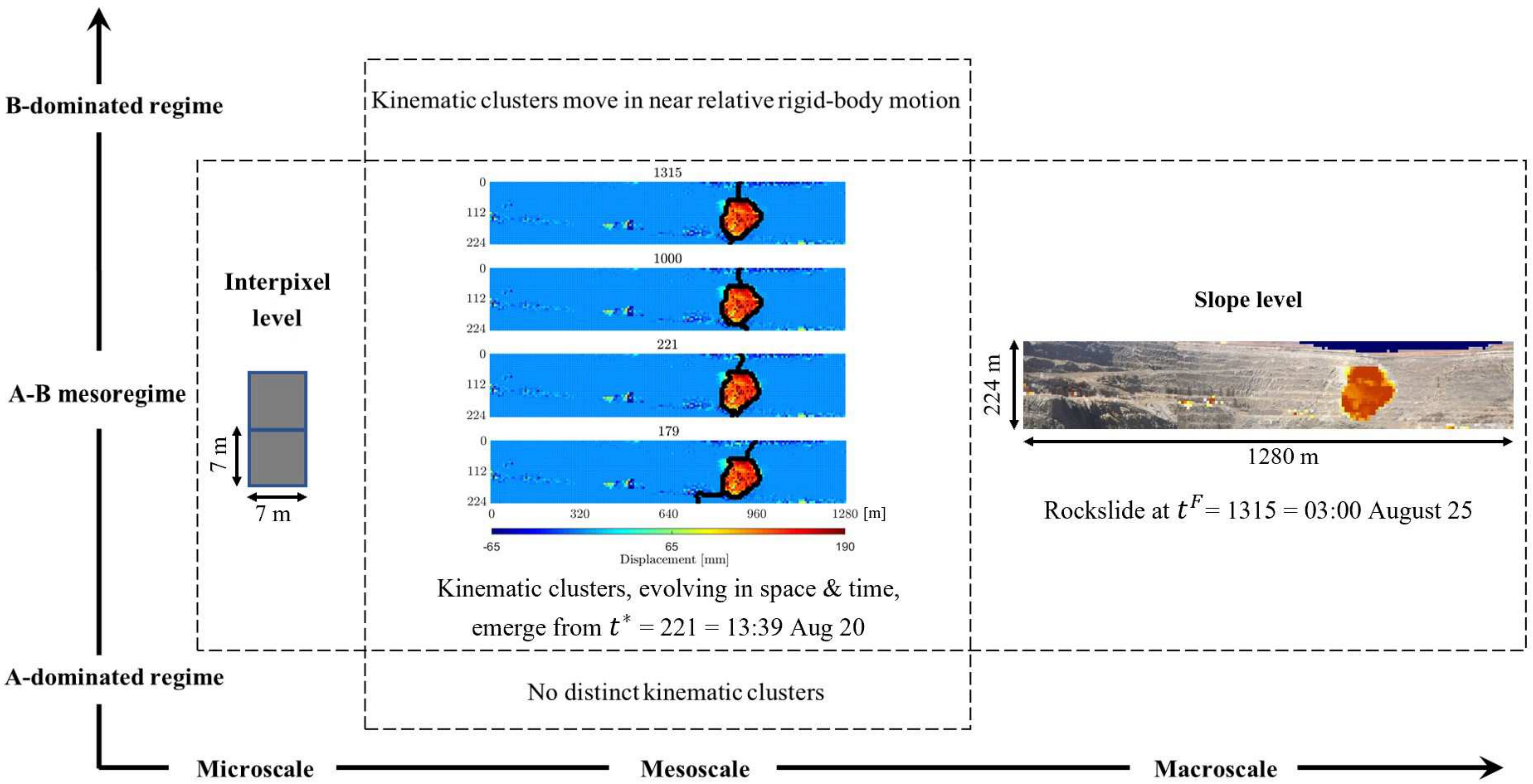}
\caption{(Color online) The precursory failure regime (PFR) over the course of monitoring a developing rockslide in an open pit mine.  This chart summarizes the mesoscience perspective of a mesoregime (PFR) where two mechanisms (A and B) coexist and give rise to emergent mesoscale kinematic clusters.  The clusters share a common boundary shown as black points overlaid on top of the displacement map at the time of the rockslide (chart-centre).  This chart is analogous to the mesoscience perspective depicted for gas- or solid- particle flow systems \cite{li2018}. }
\label{fig:mesofield}
\end{figure}

 {\it c) The principles of mesoscience. \,\,}  To account for the compromise-in-competition among preferred transmission paths and simultaneously `jump scale' -- from laboratory to field -- we integrate the network flow approach~\cite{ATSKCRMNJT,ATSKCRMNJT2} with the principles of mesoscience (Table~\ref{Tab:compete}, Figure~\ref{fig:mesofield}). Pioneered by Li and co-workers~\cite{li2013,li2018mesoscience,li2018} in the area of chemical and process engineering, mesoscience has enabled the upscaling of models of gas/solid-particle flow systems from laboratory to industrial scale.  Mesoscience is predicated on the concept of a compromise-in-competition between {\it at least} two dominant mechanisms in a so-called mesoregime of a complex system.  In the simplest case of two competing mechanisms (A and B), the mesoregime mediates two limiting regimes; each is governed by one dominant mechanism, A (B) in the A-dominated (B-dominated) regime, which is formulated as an extremum.  Li et al. argues that, while the classical single objective optimization formalism applies to each limiting regime, the compromise-in-competition in the $A-B$ mesoregime necessitates a multiobjective optimization approach. Results from prelude studies~\cite{ATSKCRMNJT,ATSKCRMNJT2}, employing a dual objective network flow analysis, corroborate this view.  
 
Moreover, opposing trends manifest as the system evolves from one limiting regime to the other ($A \rightarrow A-B \rightarrow B$ and vice versa), consistent with the mesoscience principles (Figure~\ref{fig:mesofield}, Table~\ref{Tab:compete}).  In laboratory tests where detailed analysis of underlying mechanisms are possible, the B-dominated failure regime is characterized by bursts to a peak in all the indicators of stored energy release and dissipation, including: kinetic energy, dissipation rate, population of buckling force chains and their supporting 3-cycles, average values of local nonaffine motion, grain velocity and rotation~\cite{tordesillasEtAl2016,microband,Tordesillas2007,TORDESILLAS2009706}.  By contrast, at the opposite extreme, in the A-dominated stable regime, all of these quantities are negligibly small. In PFR, these opposing tendencies compromise and give rise to spatiotemporal dynamical patterns~\cite{SDAT,KSAT,ZhouAOAD,WLF5}.

 {\it d) Clustering patterns in the kinematics characterize the mesoregime PFR. \,\,} The compromise-in-competition between force transfer and damage paths in PFR gives rise to collective motion or kinematic partitions: mesoscale clusters where constituent members move collectively in near rigid-body motion~\cite{SDAT,KSAT,ZhouAOAD,WLF5}.  
 Interestingly, Li and co-workers also observed particle clusters in the mesoregime of gas/solid-particle flow systems, and conjectured that these emerge from particles tending to minimize their potential energy, while the gas tries to choose a path of least resistance through the particle layers~\cite{li2013,li2018mesoscience,li2018}. Analogously, in the systems studied here, {\it damage favors the path of least resistance to failure -- the path that forms the common boundaries of kinematic clusters~\cite{ATSKCRMNJT,ATSKCRMNJT2}.}


 {\it e) Dynamics of kinematic clusters provide early prediction of failure across scales. \,\,} A complex network analysis of individual grain motions in sand in laboratory tests~\cite{TordesillasWalkerAndoViggiani2013} and of surface ground motion in a slope (e.g., Mine 1)~\cite{TordesillasZhouBatterham2018} has shown that the impending failure region develops in between subregions of transient but high kinematic similarity early in PFR.  Moreover, the spatiotemporal dynamics of these clusters can deliver a reliable change point $t^{*}$ from which such partitions become incised in the granular body, giving rise to their near relative rigid body motion: for example, when the active `slip region' of a slope begins to detach and accelerate downslope from a relatively stationary region below; or when parts of a rock mass on either side of a developing crack undergo relative slip. That is, persistent partitions in kinematics space forewarn of impending partitions in physical space~\cite{WLF5,KSAT,SDAT,ZhouAOAD}.  In a parallel effort~\cite{WangSS}, the computational challenges of embedding knowledge of kinematic clustering in a stochastic statistical learning model from high-dimensional, non-stationary spatiotemporal time series data were overcome, with displacement and velocity trends and the failure region of Mine1 successfully predicted more than five days in advance.

 {\it f) Establishing a connection to first principles fracture and failure mechanics for granular solids. \,\,} Relative motions at the grain-grain level were used to study the coevolution of force and damage propagation in a network flow analysis -- with explicit connections to the most popular fracture criteria, starting with Griffith's theory for crack propagation~\cite{ATSKCRMNJT,ATSKCRMNJT2}.  The emerging {\it flow bottlenecks} for force and energy, proven to be the paths of least resistance to failure, were found to deliver an accurate and early prediction of the location of shear bands and macrocracks that ultimately develop in the failure regime.  The question that now arises is: {\it Can a combined mesoscience and network flow approach detect the bottlenecks and kinematic clusters from radar-measured surface ground motion data and, if so, how can their spatiotemporal evolution be used to deliver an early prediction of a likely place and time of failure?}


Here we answer this question and demonstrate our approach through SSSAFE.  Based solely on kinematic data for input, SSSAFE first applies the network flow model to identify and characterize the emerging kinematic clusters in PFR, and then uses their dynamics to deliver an early prediction of where and when failure is likely to develop.   Different from ~\cite{TordesillasZhouBatterham2018,SDAT,KSAT,ZhouAOAD} which adopt an essentially pattern-mining approach, SSSAFE rigorously predicts the path of least resistance to failure in a manner consistent with the fundamental micromechanics and dynamics of failure across different system levels and scales. Four systems are analyzed: a standard laboratory test (Biax); and three rock slopes, man-made slopes Mine 1 and Mine 2 and a natural slope Xinmo.  The input kinematic data to SSSAFE comprise individual grain displacements in Biax, and radar line-of-sight displacement data gathered from ground-based radar (Mine 1 and Mine 2) and space-borne radar (Xinmo).

\begin{figure}
\centering
\includegraphics[width=1\textwidth]{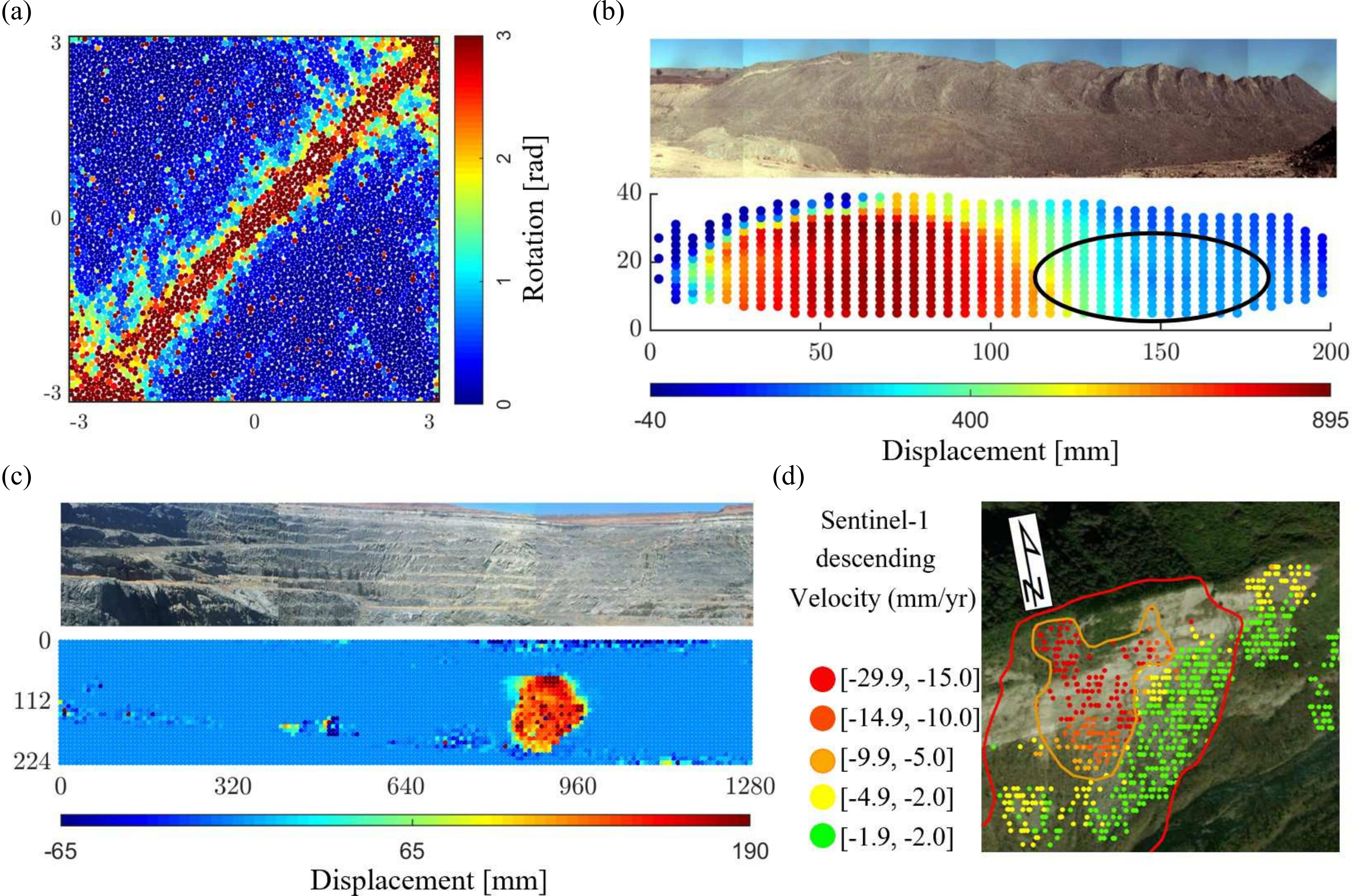}
\caption{(Color online) The systems under study in the B-dominated failure regime.  (a) Map of the cumulative absolute grain rotation in sample Biax showing the shear band where plastic deformation and energy dissipation concentrates.  Map of the cumulative line-of-sight displacement for the rock slopes, highlighting the failure location (orange-red): (b) Mine 1, (c) Mine 2 and (d) Xinmo (dimensions of Mines 1 and 2 are in meters).  In Mine 1, the second region of instability to the east (encircled) stabilized the day before the collapse.  }
\label{fig:system}
\end{figure}


\section*{Data}  
\label{sec:data}

The input data to our analysis consist of the following system properties at each time state of the monitoring period $t=1,2,....,T$: 
(a) coordinates of observation points $\ell=1,2,....,L$; (b) displacement vector recorded at each point ${\vec{d}_1,\vec{d}_2,...., \vec{d}_L}$.  
We have four data sets (Figure~\ref{fig:system}). The first is Biax, a well-studied simulation of granular failure in a standard laboratory test in which an assembly of polydisperse spherical grains is subjected to planar biaxial compression~\cite{tordesillasEtAl2016,microband,Tordesillas2007,TORDESILLAS2009706}.   Here each point $\ell$ is a moving grain and the vector $\vec{d}_\ell$ is two-dimensional.  The sample begins to dilate at around $t=50$.  Collective buckling of force chains initiate at around $t=98$, giving way to a brief period of strain-softening and the development of a single shear band along the forward diagonal of the sample.  This shear band becomes fully formed at $t=104$, referred to as the time of failure $t^F_B$ (Figure~\ref{fig:system} (a)).  From this point on, the sample exists as two clusters, in each of which constituent grains move collectively as one: two `solids' in relative rigid-body motion along their common boundary, viz. the shear band.  Details of this simulation and mechanisms underlying its bulk behavior in the lead up to and during failure are provided elsewhere~\cite{tordesillasEtAl2016,microband,Tordesillas2007,TORDESILLAS2009706}.

Three large field scale data sets are examined.  Mines 1 and 2 are from 
monitoring data of a rock slope in two different open pit mines using ground-based SSR-XT – 3D real aperture radar\cite{Dicketal2015,harries2006}  (Figure~\ref{fig:system} (b, c)). The mine operation, location and year of the rockslides are confidential.  However we have all the information needed for this analysis.  Each observation point $\ell$ is a grid cell or pixel, ranging in size from 3.5m x 3.5m to 7m x 7m, in a fixed grid.  The vector $\vec{d}_\ell$ is 1D, which corresponds to the displacement along a line-of-sight (LOS) between the radar and the point $\ell$ on the slope surface.  

Mine 1 is an unconsolidated rock slope (Figure~\ref{fig:system} (b)).  The monitored domain stretches to around 200 m in length and 40 m in height.  Movements of the rock face were monitored over a period of three weeks: 10:07 May 31 to 23:55 June 21.  Displacement at each observed location on the surface of the rock slope was recorded at every six minutes, with millimetric accuracy. This led to time series data from 1803 pixel locations at high spatial and temporal resolutions for the entire slope.  A rockslide occurred on the western side of the slope on June 15, with an arcuate back scar and a strike length of around 120 m.  Mine 1 reached peak pixel velocity of around 640 m/yr.  Considering a precautionary correction for the radar line of sight, this falls in the moderate velocity category~\cite{CrudenVarnes1996} and corresponds to an evacuation response~\cite{hungr2014}. The time of collapse $t^F_1$ occurred at around at 13:10 June 15, close to when the global average peak velocity of 33.61 mm/hr was reached. There is a competing slide: a second region of instability, to the east (encircled area, Figure~\ref{fig:system} (b)).  This region intermittently developed large movements, but the instability was somehow arrested and movement slowed down the day before the collapse of the west wall~\cite{TordesillasZhouBatterham2018,KSAT}.  In this context, this region is sometimes referred to as a false alarm in the sense that it did not eventuate into a collapse~\cite{Intrieri2017}.  While in many cases “tertiary creep” ends with a total or partial failure, it is also possible, like in Mine 1, that the whole landslide or a part of it finds a new equilibrium~\cite{intrieri2019,hungr2014}. There are many possible reasons for this, such as a reduction of the destabilizing forces through stress redistributions or the geometric configuration of the sliding surface, which slows down and ultimately arrests the whole or part of the landslide body.   There are also certain landslide types, such as earth flows~\cite{hungr2014}, which do not have a ductile behaviour and cannot experience a catastrophic failure.  Nevertheless, such flows can exhibit an exponential acceleration that, in terms of public safety, still pose an emergency challenge to deal with~\cite{intrieri2019}.

Mine 2 is a consolidated rock slope of an open cut mine dominated by intact igneous rock that is heavily structured or faulted by many naturally occurring discontinuities (Figure~\ref{fig:system} (c)). A slope stability radar scanned the section of the rock face for displacement for approximately 6 days from 15:39 August 19, until 07:05 August 25, each scan taking approximately 6 minutes, again with millimetric precision. Measurements at 5394 pixel locations were taken every 6 minutes giving high spatial and temporal resolution for the entire domain, measuring 1280 m wide and around 224 m high.  A rockslide occurred on the Southern wall at 03:00 August 25; we refer to this as the time of failure $t^F_2$ for the rest of this paper.  The area that failed, measuring approximately 135m wide and 145m high, moved over 1 million tonnes of debris.  Mine 2 reached peak pixel velocity of 2.8 m/day, which is classified as fast~\cite{CrudenVarnes1996}.

The Xinmo landslide is a rock avalanche (Figure~\ref{fig:system} (d)), composed of metamorphic sandstone intercalated with slate, that detached on June 24, 2017 and hit Xinmo village (Maoxian, China, $32^{\circ} \, 03' \, 58''$ N, $103 ^{\circ} \, 39' \, 46''$ E) causing 83 deaths and destroying 64 houses. The analyzed data set is focused only on the original source area that was located near the crest of the mountain ridge north of Xinmo village, at an altitude of 3431 m a.s.l..  As this source moved along the slope it entrained new rock material and reached an estimated volume of 13 million m$^3$ and a terminal velocity of 250 km/h~\cite{fan2017}. The site was not actively monitored at the time, but displacement data obtained from Sentinel-1 constellation, that takes periodical interferometric acquisitions of the area, have been retrospectively analyzed to determine if a forewarning would have been possible (Intrieri et al., 2018). The data used consisted in 45 SAR images in C-band (6.5 cm wavelength), at 5 m $\times$ 14 m spatial resolution, acquired along the descending orbit (incidence angle of 40.78) and spanning from 9 October 2014 to 19 June 2017 (that is five days before the failure). The pixels are of size 5m $\times$ 14 m.  Data, covering an area of 460 km$^2$, were elaborated with the SqueeSAR algorithm~\cite{ferretti2011} and comprised more than 130,000 measurement points.  Xinmo reached peak pixel velocity of around 27 mm/yr in the tertiary creep phase, which is very slow~\cite{CrudenVarnes1996} but later reached a terminal velocity of 250 km/h during failure~\cite{fan2017}.


\section*{Method}
\label{sec:methods}

\begin{figure}
\includegraphics[width=1\textwidth,clip]{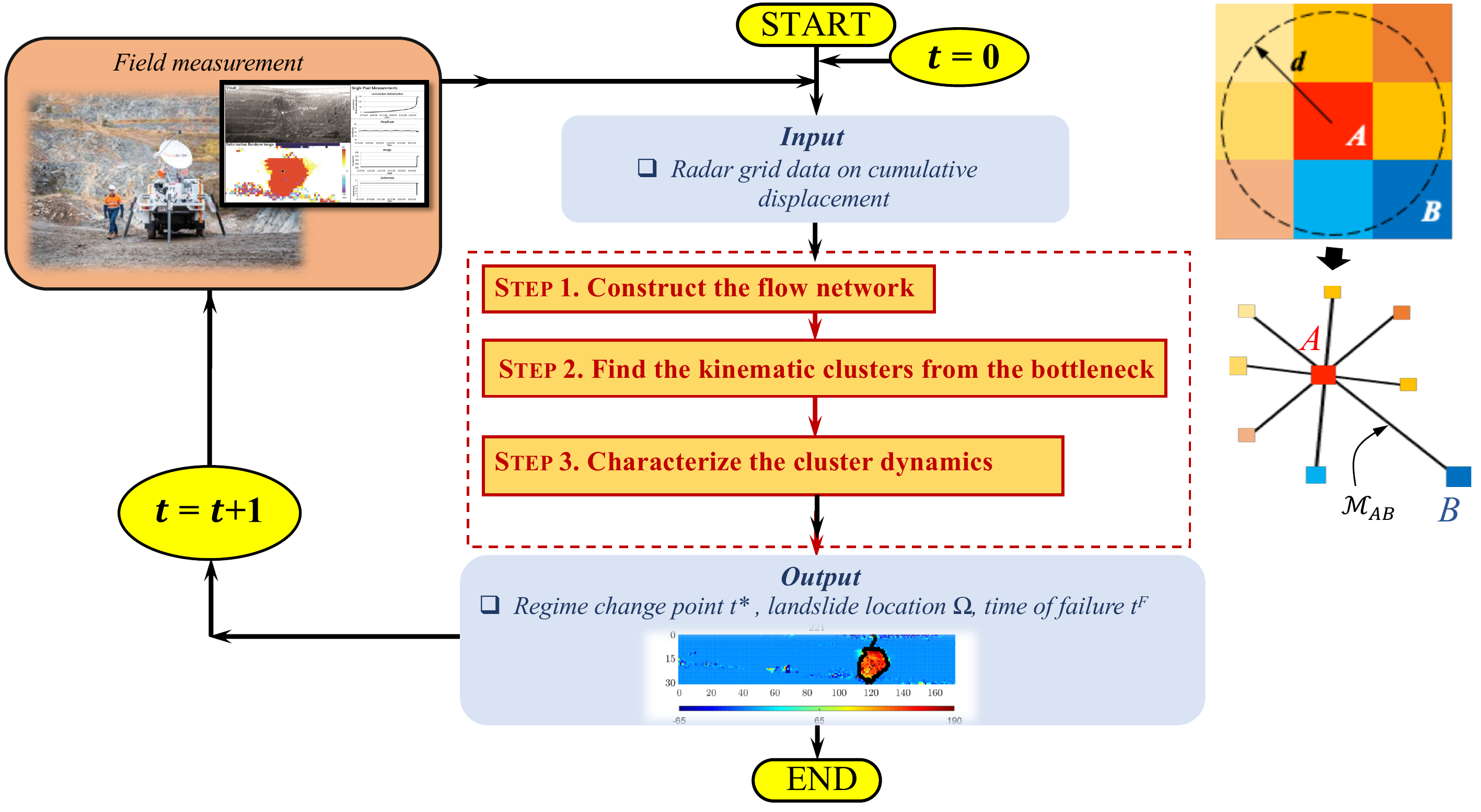}
\caption{(Color online) Flow chart summarizing the 3 steps in SSSAFE, designed for prediction of where and when failure will likely occur in a monitored domain based on spatiotemporal kinematic data.}
\label{fig:flowchart}
\end{figure}

The core components of our proposed spatiotemporal slope stability analytics for failure estimation (SSSAFE) framework are summarized in Table~\ref{Tab:scale} and Figure~\ref{fig:flowchart}.  The key idea is to model the transmission of force in each studied system in a way that accounts for the coupled evolution of the preferred pathways for force and damage, and to use this model to predict the emerging {\it kinematic clusters in the mesoregime PFR}.   To achieve a consistent formulation across different system levels, we model force transmission as {\it a flow through a network}.  At the core of this formulation is a set of optimization problems on a network in accordance with network flow theory and mesoscience principles.
We emphasize that our implementation of this model is confined only to finding the preferred paths for damage which represent the common boundary of the kinematic clusters.   Detection and characterization of the preferred paths for force are outside the scope of this investigation. Such paths have been characterized for different laboratory samples, including concrete (e.g.,~\cite{ATSKCRMNJT,tordesillas2015,PRE}).

\begin{table}
\caption{Combined mesoscience and network flow formulation behind SSSAFE.  A compromise-in-competition between mechanism A (preferred paths for force) and mechanism B (preferred paths for damage) governs the mesoscale at the laboratory level and field level.}
\centering
\begin{tabular}{cccc}
\toprule
 \textbf{System$\rightarrow$}	& \textbf{Biax}	& \textbf{Mines 1 \& 2 }   & \textbf{Flow network model}\\
 \textbf{Scale $\downarrow$} 	& 	&     &  \\
 \midrule
microscale	& grains, grain-grain interaction	& pixels, pixel-pixel interaction & nodes, links\\
mesoscale	& A versus B   & A versus B & A versus B \\
macroscale	& sample				& slope   & flow network\\
 \bottomrule
 \end{tabular}
 \label{Tab:scale}
 \end{table}

\subsection*{Core components of SSSAFE}
\label{secFlow}

The core components of SSSAFE are implemented in three consecutive steps, with Steps 2 and 3 delivering respectively the predictions on the likely location and the time of failure.  Recent work~\cite{ATSKCRMNJT2} shows explicit connections between the formulation below and the most popular fracture criteria of fracture mechanics, beginning with Griffith's theory for crack propagation.  The method below, consistent with Griffith's theory, was found to provide the most accurate and early prediction of failure in PFR for a range of ductile to quasi-brittle laboratory samples.

{\bf Step 1: Construct the flow network $\mathcal{F}$.}  Forces are transmitted along physical connections.  Hence the construction of the flow network $\mathcal{F}$ begins with an undirected network $\mathcal{N}$ that represents the physical connectivity of the system: the grain contact network in Biax or the proximity network in Mines 1 and 2 where pixels within a distance $d$ of each other are connected.  Each node of $\mathcal{N}$ represents a grain (pixel), while each link in $\mathcal{N}$ represents a grain-grain contact (pixel-pixel connection).  The links in $\mathcal{N}$ vary with loading history in Biax, but is fixed across the monitoring period for both Mines 1 and 2.

\begin{figure}
\centering
\includegraphics[width=0.7\textwidth,clip]{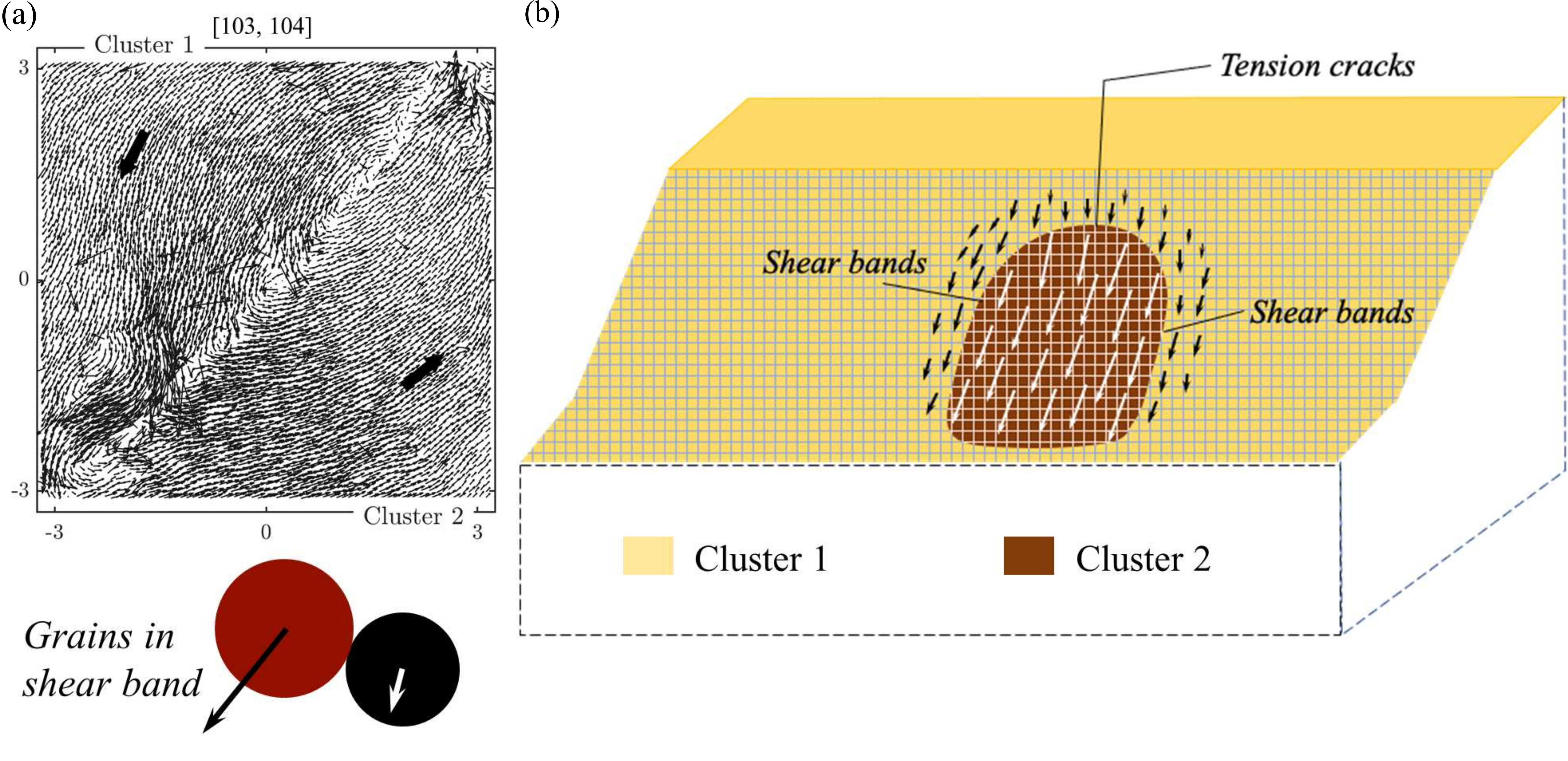}
\caption{(Color online) Collective motion of mesoscale clusters characterizes the terminal stages of the mesoregime PFR. Emerging kinematic clusters increasingly move in near relative rigid-body motion: (a) the actual displacement field at failure in Biax, (b) depiction of surface ground motion on a slope. Links along the shared boundary of kinematic clusters, $\Omega$, are closest to breaking point (i.e., smallest total path capacity $c(\Gamma)$) due to the large relative motions of its constituent elements. }
\label{fig:pattern}
\end{figure}

Next, $\mathcal{N}$ is transformed to a directed network $G = (V, A)$, where $V$, $A$ are the set of nodes and set of arcs respectively.   That is, each link connecting nodes $i \in V$ and $j \in V$ in $\mathcal{N}$ is represented by a pair of symmetric arcs $e \in A$: one from $i$ to $j$ and another vice versa.  Given this symmetry, we will use the symbol $e$ to also denote a link.  Every link in $G$ is then assigned a non-negative {capacity} $c(e)$ which corresponds to the maximum flow value that the given link can support.  Since the model concerns force transmission, $c$ is thus the force that must be overcome to break the connection: the strength or {\it resistance to failure} of the grain-grain contact or pixel-pixel connection.  Given that what is measured reliably in both laboratory and field levels is motion -- not forces or stresses -- we express this capacity $c$ in terms of the motions of the connected elements.  Hence, the contact capacity function $c$ is given by 
\begin{equation}
    c(e)=c_{ij} =c_{ji} = \frac{1}{|\overrightarrow{\Delta u_{ij}}|^2},
    \label{eq:cap}
\end{equation}
where $|\overrightarrow{\Delta u_{ij}}|$ is the magnitude of the relative displacement of two grains (or two pixels) linked in $\mathcal{N}$.  
Note that since we are only interested in the flow bottleneck~\cite{ATSKCRMNJT,AhujaNetworkFlows}, what is important in this analysis are the relative values of the link capacities and not their absolute values.  Indeed, for this purpose, the model for the link capacity need not be in units of force, as previously shown (e.g.,~\cite{ATSKCRMNJT,tordesillas2015, PRE}). Consequently, in Equation~\eqref{eq:cap}, we set the capacity to be such that the higher the relative motion of grains (pixels) linked in $\mathcal{N}$, the less stable is the connection and in turn the lower is its corresponding capacity $c$ (see Figure~\ref{fig:pattern}).  Finally, a direction for the flow is dictated by a pair of artificial nodes $q$ and $k$ called the \emph{source} and the \emph{sink} of $G$.  The quadruple $\mathcal{F} = (G, c, q, k)$ is called a {\it flow network}.  In the Biax sample, the natural choice for the \emph{source} $q$ and \emph{sink} $k$ are the top and bottom walls so that the direction of flow is in alignment with the direction of the applied vertical compression (and major principal stress axis) of the sample.


{\bf Step 2: Find the kinematic clusters from the bottleneck of $\mathcal{F}$.}   
The bottleneck of $\mathcal{F}$, $B(\mathcal{F})$, is given by the cut of $\mathcal{F}$ with the least capacity.  Any cut of $\mathcal{F}$, $\Gamma$, is a set of links in $\mathcal{N}$ which, if disconnected, represents a literal cut of $\mathcal{F}$ into two disjoint components $\{W,W'\}$ of $V$ such that no flow can be transmitted from source $q\in W$ to sink $k\in W'$.  Thus, any cut $\Gamma$ contains all arcs emanating from a node in $W$ and terminating on a node in $W'$.  

Physically, a cut $\Gamma$ may be thought of as a virtual crack of the studied granular body or domain whose connectivity is described by $\mathcal{N}$.  Physical disconnection of the contacts associated with the links in $\Gamma$ would thus result in a literal system-spanning crack which splits the body into two disjoint pieces.  The capacity of $\Gamma$ is defined as $c(\Gamma) = \displaystyle \sum_{e\in\Gamma}{c(e)}.$  Following Equation~\eqref{eq:cap}, this represents the total force flow that must be overcome to disconnect every link in $\Gamma$.  

Here we are interested in finding that cut with the least capacity -- the so-called minimum cut, also known as the bottleneck $B(\mathcal{F})$.  Thus, the capacity of the bottleneck $B(\mathcal{F})$ represents the global failure resistance, $F^*$: the minimum amount of force flow needed to overcome the resistance of the connected links $B(\mathcal{F})$ to break apart and split the granular body into two disjoint pieces.  Note that this analysis does not preclude a body from splitting apart into more than two pieces: in such cases, one can repeat the same analysis described here for each piece to obtain further subpartitions.  In the cases studied here, this is unnecessary as the studied systems split apart essentially into two components with the bottleneck being their shared boundary.   
In Biax, the bottleneck $B(\mathcal{F})$ predicts the location of the shear band that forms in the failure regime.  In the case of Mines 1 and 2 and Xinmo, $B(\mathcal{F})$ predicts the boundary of the landslide.  As time to failure draws near, we expect motion in the components to become increasingly coherent and near-rigid-body resulting in kinematic clustering.  The {\it active cluster} in PFR, denoted by $\Omega$, distinguishes itself by manifesting an increasing downward motion (viz. increasing trend in cumulative displacement and velocity) due to gravity, while the stable cluster remains relatively stationary.

To find the bottleneck of Biax at each time, we solve the \textsc{Maximum flow - Minimum cut (MFMC) problem} on $\mathcal{F}$, following earlier work~\cite{ATSKCRMNJT,PRE}.  This is a two stage calculation.  Stage 1 solves the \textsc{ Maximum flow problem} to find the global flow capacity, $F^*$, the maximum flow that can be transmitted through $\mathcal{N}$ given its topology and link capacities.  More formally, given a flow network $\mathcal{F} = (G,c,q,k)$, a link flow $x(e)$ is called a {\it feasible} $q$-$k$ flow, if it satisfies:

\noindent
(a) the conservation of flow
\begin{equation}
\sum_{e \in \d^{-}(v)} x(e) \;= \;\sum_{e \in \d^{+}(v)} x(e), \; \; \; \forall v \in V - \{q, k\}, 
\label{eq:flow-con} 
\end{equation}
where ${e \in \d^{-}(v)}$ denotes arcs entering node $v$ and ${e \in \d^{+}(v)}$ denotes arcs leaving node $v$;

\noindent
(b) the capacity rule
\begin{equation}
0 \le x(e) \le c(e),\; \; \; \; \forall e \in G.
\label{eq:flow-caps}
\end{equation}
Hence, the force flow along each link, $x(e)$, is regulated by the threshold for damage $c(e)$ which is a function of the relative motion between the connected elements (Equation~\ref{eq:cap}).

\begin{enumerate}
\item[]The \textsc{Maximum Flow Problem} can be expressed as: find a feasible $q$-$k$ flow $x$ such that the following flow function $f(x)$ is maximum:
\begin{equation}
f(x) = \sum_{e \in \d^{+}(q)} x(e) - \sum_{e \in \d^{-}(q)} x(e).
\label{eq:flow-vals} 
\end{equation}
The {\it flow value} that solves the above is the maximum flow $F^*$. 
\end{enumerate}
\noindent 
Once $F^*$ is established, we move to Stage 2 to solve the \textsc{Minimum Cut Problem}.  

\begin{enumerate}
\item[]The \textsc{Minimum Cut Problem} of $\mathcal{F}=(G,c,q,k)$ is the cut $\Gamma_{min}$ such that
\begin{equation}
c(\Gamma_{min}) = Minimize\big\{ \displaystyle \sum_{e\in\Gamma}{c(e)} \big\}.
\label{eq:capacity} 
\end{equation}
\end{enumerate}
\noindent

The above is typically solved using the Ford-Fulkerson algorithm \cite{AhujaNetworkFlows}.  This exploits the well known \emph{max-flow min-cut theorem} which states that the maximum flow possible $F^*$ is the capacity of the minimum cut or bottleneck~\cite{liu2011segmenting}.  
Using this theorem and Equations~\eqref{eq:cap} -- \eqref{eq:capacity}, we can now directly relate the conditions on where and when catastrophic failure occurs to the bottleneck.
That is, catastrophic failure occurs when the force flow exceeds the resistance to breakage of all the links in the bottleneck.  Where the system physically breaks apart is given by the bottleneck itself, which is the limiting shared boundary of the kinematic clusters.

\begin{figure}
\centering
\includegraphics[width=.6\textwidth,clip]{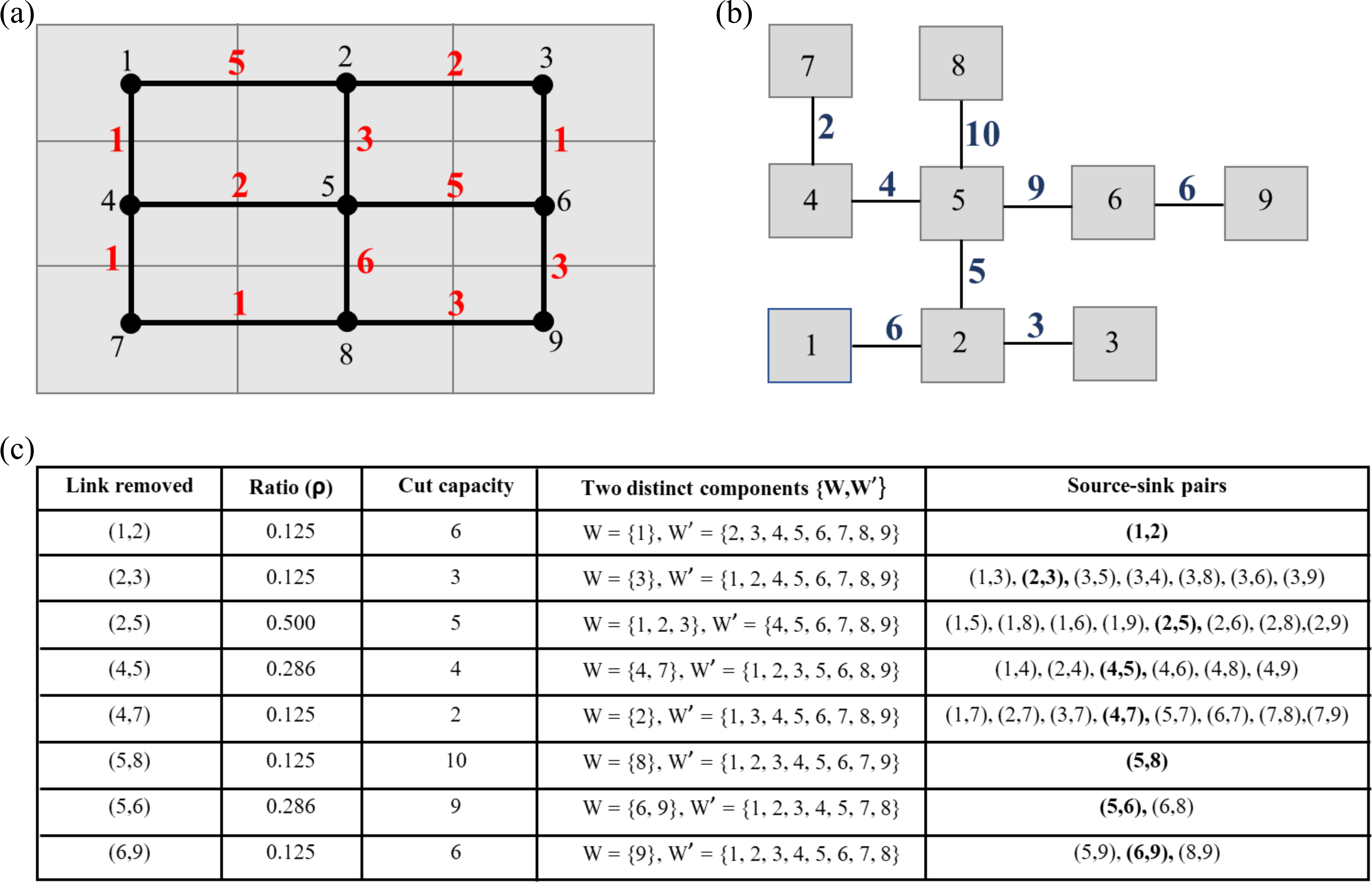}
\caption{(Color online) (a) An example contact network $\mathcal{N}$ of $n = 9$ nodes (pixels) with contact capacities shown in red.  (b) Corresponding Gomory-Hu tree $\mathcal{T}$ for $\mathcal{N}$ in (a). (c) Summary of the explicit source-sink pairs in $\mathcal{T}$ with 8 distinct minimum cuts and their related properties including the resulting clusters or components in $\mathcal{N}$ when the link of each pair is removed in $\mathcal{T}$.  The last column lists all implicit source-sink pairs with the same minimum cut as the explicit source-sink pair in bold.}
\label{fig:stM1}
\end{figure}

While finding the bottleneck in Biax is relatively straightforward, this is not the case for Mines 1 and 2.   The difficulty arises because there is no obvious choice for the source-sink pair to direct the flow, given the uncontrolled and unknown loading conditions of these slopes.  To address this, we construct the Gomory-Hu tree (GHT)~\cite{gomory1961multi} for the network $\mathcal{N}$.  The procedure is described in \cite{KahagalagePhDthesis} but in what follows we outline this briefly for completeness.   Let $G^* = (\mathcal{N} ,c)$ be an undirected, link-capacitated network.  For every pair of nodes $u,v$ in $G^*$, the GHT stores information on the minimum $u$-$v$ cut of $G^*$ that separates $u$ and $v$.  For a network with $n$ nodes, there are a total of $n(n-1)/2$ possible source-sink pairs, each with a corresponding minimum cut. However, the construction of the GHT shows that the minimum cuts for some pairs of nodes are identical. In fact, the GHT contains information on exactly $n-1$ distinct minimum cuts~\cite{gomory1961multi} corresponding to a set of $n-1$ explicit source-sink pairs.  One could infer all the remaining implicit source-sink pairs from this set using the GHT, as illustrated in Figure~\ref{fig:stM1}.  Formally, the GHT is defined as follows. 


\begin{definition}
\label{def:ght}
\textit{For a given link-capacitated network, $G^* = (\mathcal{N},c)$, a tree $\mathcal{T}$ is a Gomory-Hu tree (GHT) if the following holds:}
\begin{enumerate}
\item The nodes of $\mathcal{T}$ coincide with the nodes of $G^*$.
\item Each link $l$ in $\mathcal{T}$ has a non-negative weight $w(l)$.
\item For each pair of nodes $u,v$ in $\mathcal{T}$, let $l_{m}$ be the link of minimum weight on the path joining $u$ and $v$ in $\mathcal{T}$.
\end{enumerate}
Then $w(l_{m})$ is equal to the capacity of the minimum cut separating $u$ and $v$ in $G^*$.
\end{definition}

In Figure~\ref{fig:stM1}, we illustrate an example contact network $\mathcal{N}$ with $n = 9$ pixels, its corresponding Gomory-Hu tree $\mathcal{T}$ and a table summarizing the outcome of removing a link in $\mathcal{T}$. There are 36 possible source-sink pairs. $\mathcal{T}$ contains 8 explicit source-sink pairs (column 1, Figure~\ref{fig:stM1} (c)). Removing link $l$, connecting nodes $u$ and $v$ in $\mathcal{T}$, gives two distinct components $\left\lbrace W, W' \right\rbrace$: these correspond to the kinematic clusters of $\mathcal{N}$ when the edges in the minimum cut separating the source-sink pair $u$ and $v$ are removed.  All other source-sink pairs and their corresponding minimum cuts can be inferred from $\mathcal{T}$.  

Consider, for example, the minimum cut of $\mathcal{N}$ separating the source-sink pair $u=1$ and $v=8$.  Link $l_m=(2,5)$ has the minimum weight in the path from $u=1$ to $v=8$ in $\mathcal{T}$ (Definition~\ref{def:ght}).  Thus, removing $l_m=(2,5)$ $\mathcal{T}$ results in $W=\left\lbrace1,\;2,\;3\right\rbrace$ and $W'=\left\lbrace4,\;5,\;6,\;7,\;8,\;9\right\rbrace$.  In $\mathcal{N}$, this partition corresponds to the removal of edges $(1,4)$, $(2,5)$, and $(3,6)$ that constitute the minimum cut for the source-sink pair $u=1$ and $v=8$ with capacity of 5.  Note that there are other source-sink pairs having the same minimum cut.  

From $\mathcal{T}$, the absolute (global) minimum cut capacity is 2. The corresponding two partitions are $W = \left\lbrace7\right\rbrace$ and $W' = \left\lbrace1,\; 2,\; 3,\; 4,\; 5,\; 6,\; 8,\; 9\right\rbrace$. The global minimum cut contains edges $\left\lbrace(4,7),\; (7,8)\right\rbrace$.  Observe this global minimum cut is biased towards highly imbalanced cuts where one component is significantly smaller than the other in terms of the number of member nodes.  Such highly imbalanced partitions may correspond to the smaller component having only one or at most a few pixel locations out of thousands or more.  As such, this may not provide a complete summary of emerging partitions that lead to catastrophic failure.  Larger partitions that span the system, where the part that dislodges from the rest of the slope constitutes a sizeable portion of the slope, are of interest.  Accordingly, we introduce a cut ratio $\rho$ which is defined as the number of nodes in the smallest to largest component upon removal of a link in $\mathcal{T}$.  Hence, in this example, if we are interested in the smallest component containing at least 3 pixels, we find the minimum cut such that $0.3 \leq \rho \leq 1$.  This yields the cut that corresponds to the removal of link $(2,5)$ in $\mathcal{T}$; the explicit source-sink pair is $(u=2, v=5)$ as before. By inspecting $\mathcal{T}$, we can see that source-sink pairs $(u=1, v=5), (u=1, v=8), (u=1, v=6), (u=1, v=9), (u=2, v = 8), (u=2, v = 6), (u=2, v = 9)$ correspond to the other minimum cuts also satisfying $0.3 \leq \rho \leq 1$. Note that this requires enumeration of all possible source-sink pairs and their minimum cuts.  An outline of this procedure is given in Algorithm~\ref{alg:bcut}.  


\begin{algorithm}
\caption{The major crack prediction algorithm}
\label{alg:bcut}
\begin{algorithmic}[1]
     \STATE Let $\mathcal{T}$ be a GHT constructed for a capacitated graph $G^*$ with $n$ nodes
     \STATE Set $O$ to an $(n-1) \times 3$ matrix whose entries are all zeros 
     \FOR{$i := 1:n-1$}
         \STATE Remove a link $l_i$ in $\mathcal{T}$
          \STATE Compute the ratio $\rho(l_{i})$ upon the removal of the link $l_i$ and set the $i^{th}$ row of $O$, $O[i,:] := [\rho(l_{i}) \;w(l_i) \;l_i]$ where $w(l_i)$ is the weight of the link $l_i$ in $\mathcal{T}$
        \ENDFOR
    \RETURN The edge $l_i$ which corresponds to the minimum capacity $w(l_i)$ for $\rho_m \leq \rho(l_i) \leq 1$. Two nodes in $l_i$ gives the location of the source and the sink node and removal of this link $l_i$ in $\mathcal{T}$ gives the corresponding two components $\left\lbrace W, W' \right\rbrace$. 
\end{algorithmic}
\end{algorithm}

\begin{figure}
\centering
\includegraphics[width=1\textwidth,clip]{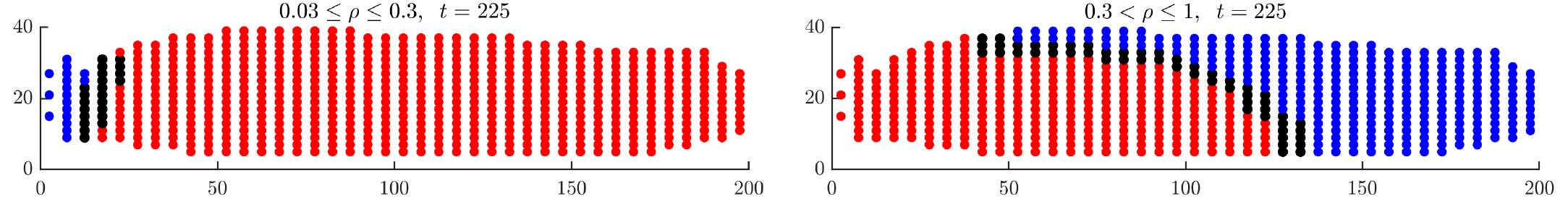}
\caption{(Color online) Kinematic clusters (red and blue) for Mine 1 at stage $t=225$ for $\rho_m \leq \rho \leq 1$.  The active cluster $\Omega$ is colored red.  Black points highlight the pixels connected by the set of links in the bottleneck, the common boundary of the clusters.  }
\label{fig:clusters}
\end{figure}

For Mines 1 and 2 and Xinmo, the number of nodes in $G^*$ and possible source-sink pairs are, respectively: $612$, $186,966$; $5394$, 14,544,921; $610$, 185,745.    It is thus computationally expensive to enumerate every source-sink pair. However, we can use $\mathcal{T}$ to capture a failure event such that the minimum cut identifies a failure area that is no smaller than a prescribed fraction of the studied domain size.  
That is, we can remove the link in $\mathcal{T}$ with the minimum weight such that $0.3 < \rho \leq 1$ to obtain two corresponding cluster components $\left\lbrace W, W' \right\rbrace$.   It is easy enough to check the partitions, and ensure that parts of landslide boundaries that are close to the boundary of the monitored domain (recall the top left boundary of the rockfall in Mine 1 (Figure~\ref{fig:system} (b))) are also captured.  This is done simply by checking the cases where $0.03 \leq \rho \leq 0.3$, which identify partitions that lead to the smaller cluster being as small as 3\% of the number of nodes in the larger cluster (Figure~\ref{fig:clusters}).  In summary, the key output from Step 2 is the bottleneck $B(\mathcal{F})$, the path of least resistance to failure which separates the active cluster $\Omega$ that will likely collapse from the rest of the slope, and the slope's failure resistance $F^*$. 

{\bf Step 3: Characterize the cluster dynamics. } As depicted in Figure~\ref{fig:flowchart}, at each time state up until the current time $t$, we find the flow bottleneck $B(\mathcal{F})$, its associated clusters and the failure resistance.  We can use this historical information to characterize the dynamics of the cluster motions as the monitoring advances in time.  Here we are interested in one of the defining aspects of granular failure, namely, collective motion. As time advances towards the failure regime, we quantify the extent to which: (a) intracluster motions become increasingly coherent and similar -- at the same time as intercluster motions become more and more different (separated in kinematic state space); and (b) the predicted clusters no longer change in member elements, suggesting that the pattern of impending failure has become physically incised in the system.  To do this, we compute the silhouette score $S$ \cite{Rousseeuw1987} to quantify the quality of the clustering pattern obtained from the network flow analysis, coupled with an information-theoretic measure of Normalized Mutual Information (NMI) \cite{Vinh10} to quantify the temporal persistence of the clustering pattern.  

\begin{figure}
\centering
\includegraphics[width=0.6\textwidth,clip]{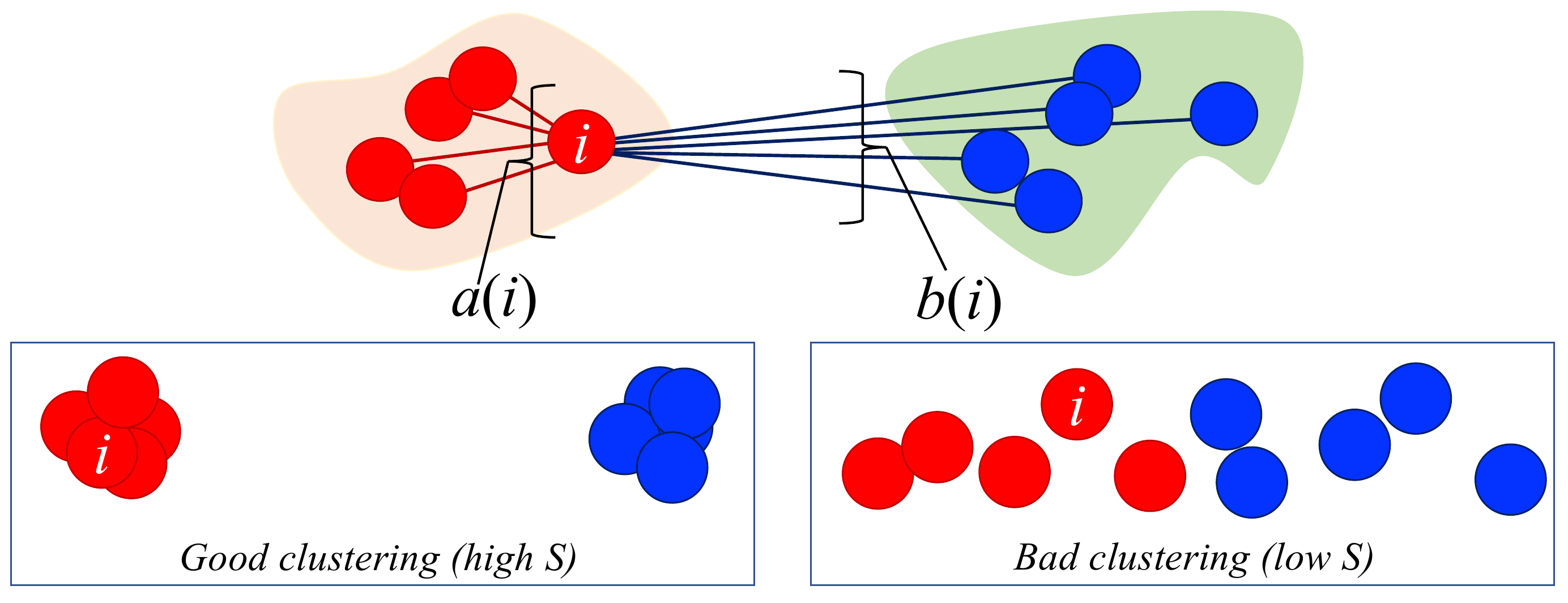}
\caption{(Color online) Depiction of the silhouette score $s(i)$ for node $i$, used to quantify the quality of clustering in kinematic state space. $a(i)$ ($b(i)$) measures intra- (inter-) cluster similarity of node $i$. }  
\label{fig:sil}
\end{figure}

The Silhouette score $S \in [-1,1]$ gives an overall measure of the quality of clustering \cite{Rousseeuw1987}.  It is the global average of $s(i)$ which measures how similar is a given node $i$ to the other nodes $j$ in its own cluster (cohesion) compared to the nodes in the other clusters (separation):
 \begin{equation}
S= {\frac{1}{n}}\sum\limits_{i=1}^n s(i) =\sum\limits_{i=1}^n \frac{b(i)-a(i)}{\text{max}\big[b(i),a(i)\big]};
 \label{eqn:Silh}
 \end{equation}
here $a(i)$ is the average distance in the displacement state-space from $i$ to all other nodes in the same cluster, and $b(i)$ is the average of the distances from $i$ to all points in the other cluster.  As shown in Figure \ref{fig:sil}, a good clustering pattern (high $S$) is one where the nodes in the same cluster exhibit very similar features (nodes are tightly packed in feature state space hence small $a(i)$); while nodes from different clusters have very different features (red nodes are well separated from blue nodes in feature state space hence large $b(i)$).  As a general guide, values below 0.2 suggest essentially no clustering pattern was found, while the closer $S$ is to one, the more compact are the individual clusters while being more separated from each other.   Given the studied feature is motion, an increasing trend with respect to time in $S$ from around 0.2 to its upper bound of 1 suggests that the clusters are moving in increasingly relative rigid-body motion.  

The Normalized Mutual Information (NMI) \cite{Vinh10} basically tells us how much knowing the clustering pattern at the previous time, $X(t-1)$, reduces our uncertainty of the clustering at the current time, $X(t)$.  The Normalized Mutual Information (NMI) is defined as
 \begin{equation}
\text{NMI}= \frac{I(X(t);X(t-1))}{\sqrt{(H(X(t))H(X(t-1)))}};
\label{eqn:nmi}
\end{equation}     
here $I(X(t);X(t-1))$ is the mutual information between $X(t)$ and $X(t-1)$  and $H(.)$ is the entropy of the corresponding clustering assignments.  NMI $\in [0, 1]$: 0 means there is no mutual information, as opposed to 1 where there is perfect correlation or similarity, between the clusters at $t$ and $t-1$.
Intuitively, NMI measures the information that the clustering assignments $X(t)$ and $X(t-1)$ share: the higher the NMI, the more useful information on the clustering pattern is encoded in $X(t-1)$ that can help us predict the clustering at the next time state $X(t)$.  


In summary, based on the results from Steps 2 and 3, we can identify a regime change point $t^*$ from which the failure resistance $F^*$ drops close to its minimum of zero, as $S$ rises and/or levels above 0.2, while NMI stays close to 1.   For all $t \geq t^*$, a prediction on the landslide region is given by $\Omega $, the active or fastest moving cluster.  In addition, the time of failure $t^F$ can be predicted by performing a linear regression with a rolling time window of the inverse mean velocity of $\Omega$ for $t \geq t^*$.  This not only obviates the need to subjectively select a pixel to implement the Fukuzono INV analysis~\cite{fukuzono1985} but also ensures this analysis takes into account the spatiotemporal and coupled evolution of force and damage pathways in PFR.

\section*{Results and discussion}
\label{sec:results}

\begin{figure}
\centering
\includegraphics[width=\textwidth,clip]{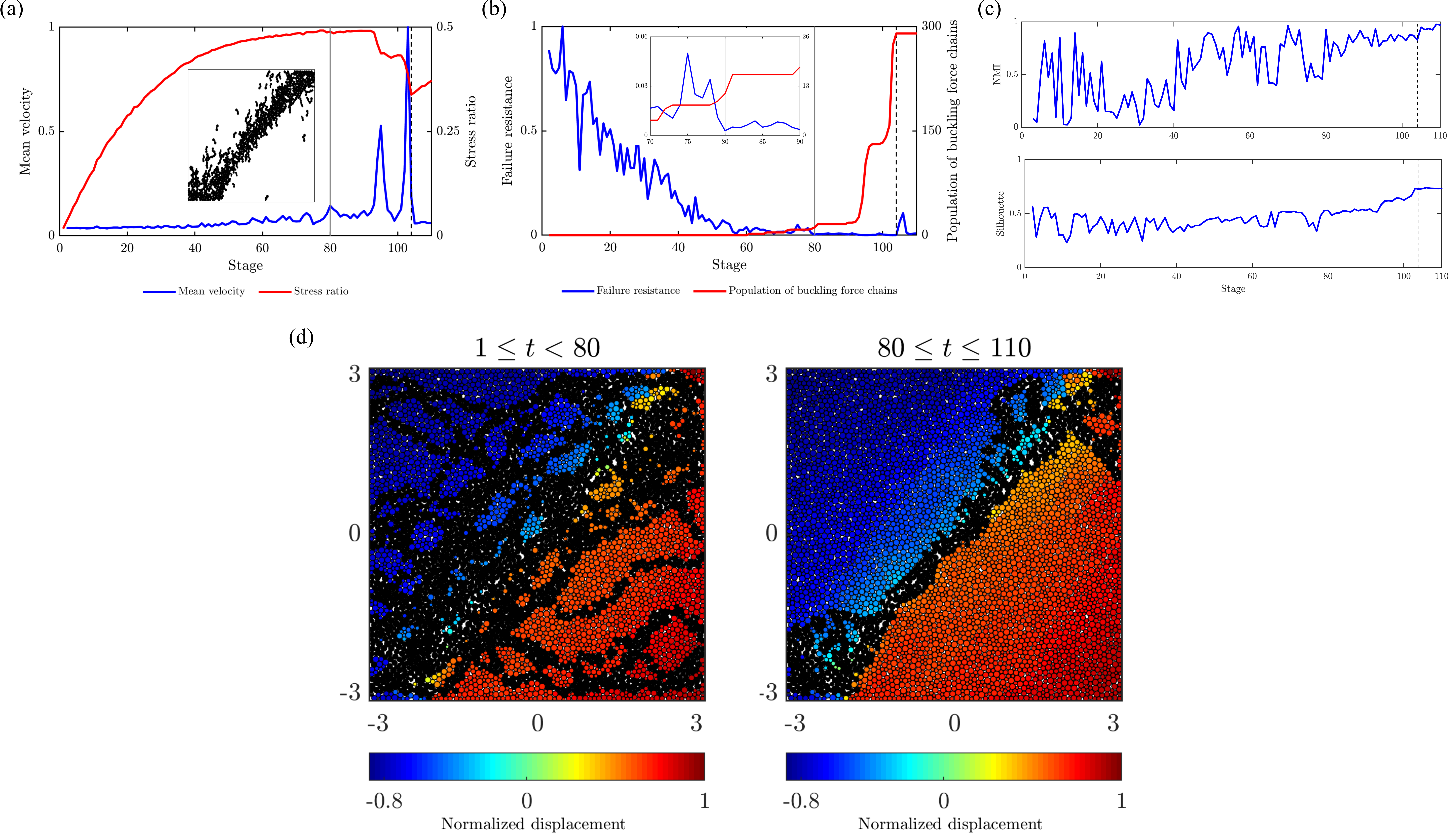}
\caption{(Color online) The mesoregime mediates the stable regime and the failure regime in Biax.  Time evolution of: (a) the Biax mean velocity, along with the shear strength of Biax as measured by the stress ratio; (b) the failure resistance $F^*$ from Step 2.   Inset in: (a) shows the collective buckling of force chains in the shear band; (b) shows zoomed-in area near the regime change point $t^*_B$= 80. (c) Plot of the time evolution of NMI and $S$ from Step 3.  Vertical lines mark the regime change point $t^*_B$=80 (solid grey line) and the time of failure $t^F_B$= 104 (dashed black line), respectively.  (d) Cumulative predictions of the shear band (black points) from Step 2 overlaid on top of the map of the magnitude of displacement at the time of failure.  Southwesterly (northeasterly) displacement is given a negative (positive) sign.   }
\label{fig:oppose}
\end{figure}

\begin{figure}
\centering
\includegraphics[width=\textwidth,clip]{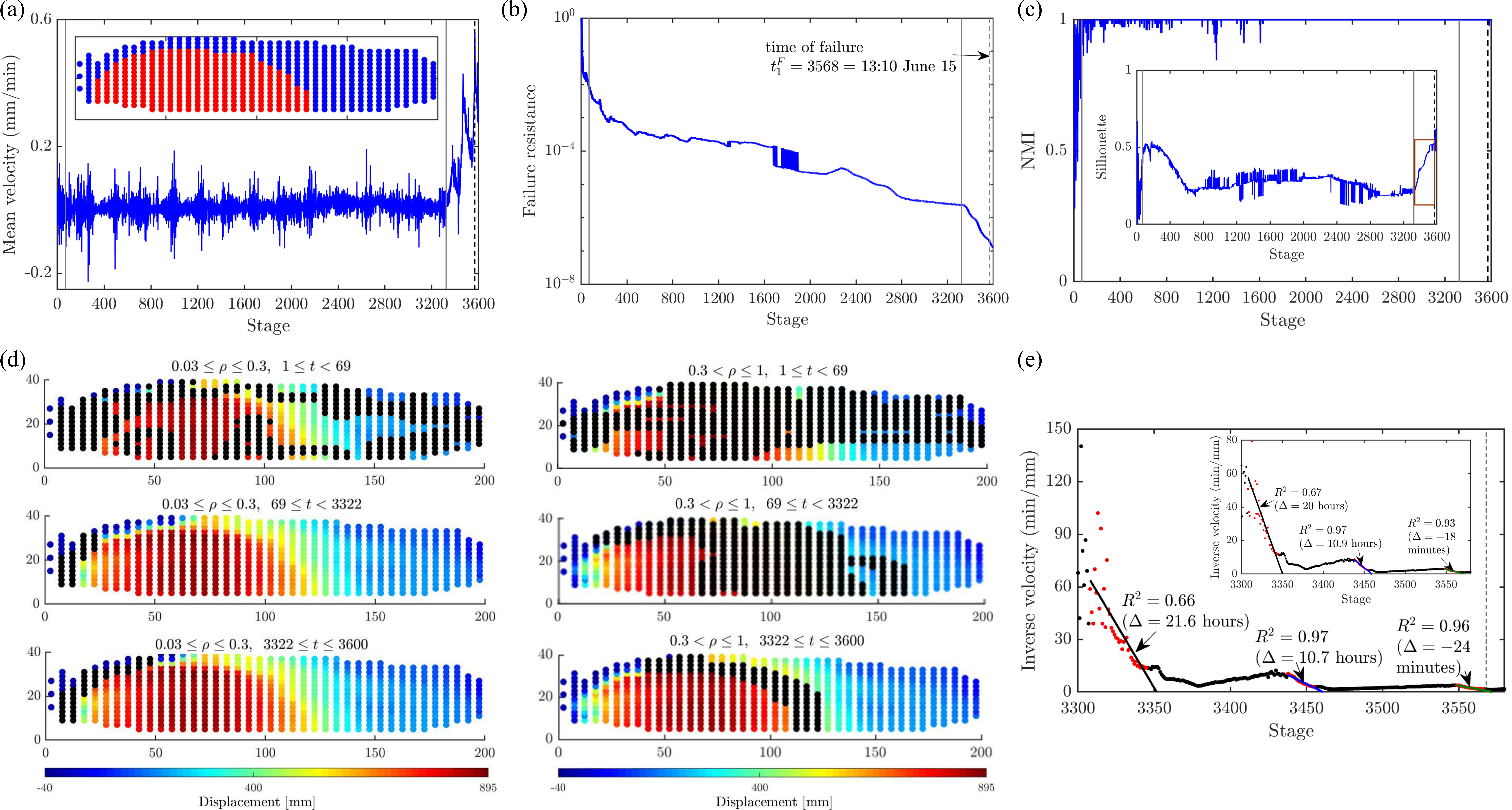}
\caption{(Color online) Slope Mine 1. Time evolution of: (a) the mean velocity of Mine 1 with the failure location in red (inset); (b) the failure resistance $F^*$ on log-axis from Step 2; (c) the NMI index and $S$ (inset) from Step 3; (e) the inverse mean velocity of $\Omega$ and of pixel $p$ (inset), over the time interval indicated by the red window in the inset in (c). Vertical lines mark the regime change points (solid grey line) $t^*_{1a}=69=$ 17:14 May 31 and $t^*_{1b}=3322=$ 12:14 June 14, and the time of failure (dashed black line) $t^F_1=3568=$ 13:10 June 15. (d) Cumulative predictions of the preferred paths for damage, and the landslide boundary for $t \geq t^*$, are highlighted by the black points and overlaid on top of the displacement map at $t^F_1$. Supplementary Movie Mine1 shows the evolution of $\Omega$ over the period of the monitoring campaign.}
\label{fig:mine1}
\end{figure}

\begin{figure}
\centering
\includegraphics[width=\textwidth,clip]{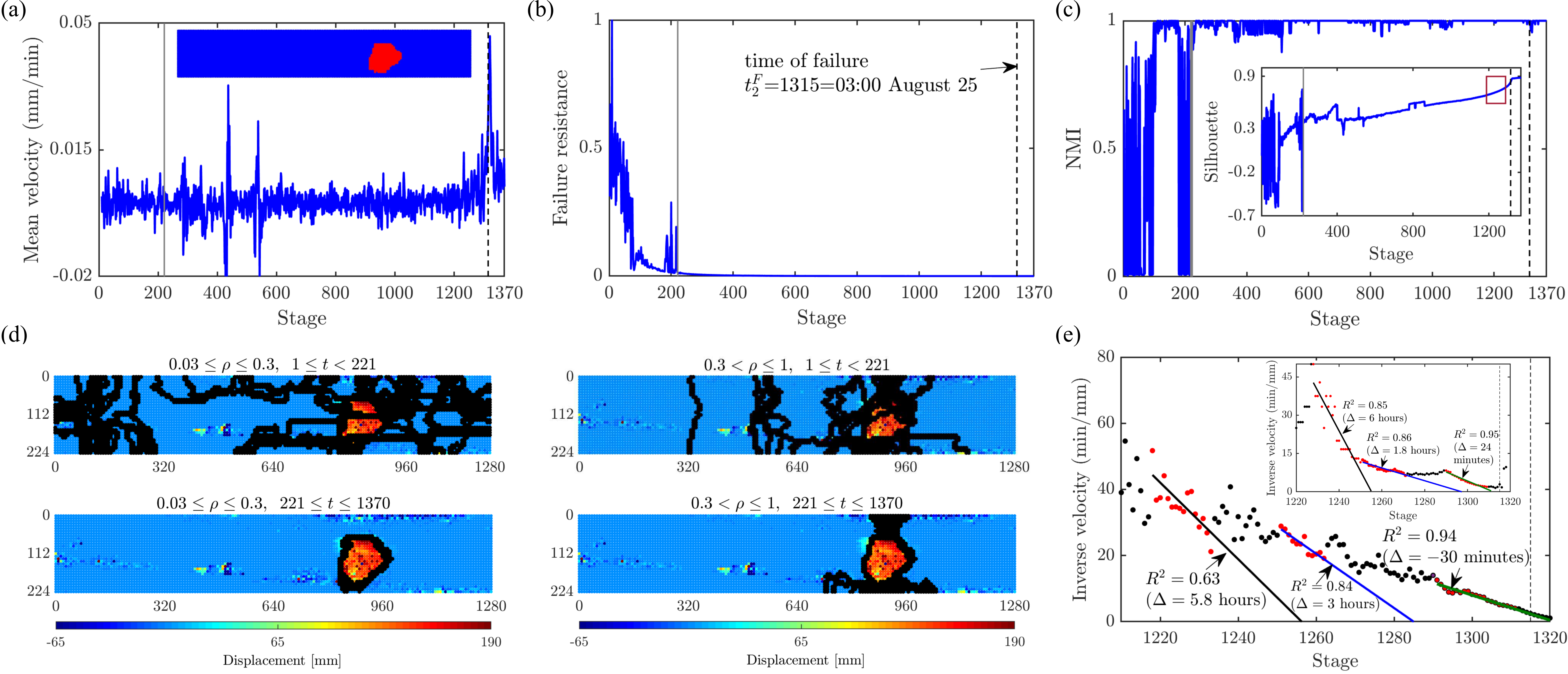}
\caption{(Color online) Slope Mine 2. Time evolution of: (a) the mean velocity of Mine 2 with the failure location in red (inset); (b) the failure resistance $F^*$ from Step 2; (c) the NMI index and S (inset); (e) the inverse mean velocity of $\Omega$ and of pixel $p$ (inset), over the time interval indicated by the red window in the inset in (c).  Vertical lines mark the regime change point $t^*_2$=221=13:39 August 20 (solid grey line) and time of failure $t^F_2$=1315=03:00 August 25 (dashed black line).  (d) Cumulative predictions of the landslide boundary (black points) are overlaid on top of the displacement map at $t^F_2$.  Supplementary Movie Mine2 shows the evolution of $\Omega$ over the period of the monitoring campaign.}
\label{fig:mine2}
\end{figure}

\begin{figure}
\centering
\includegraphics[width=\textwidth,clip]{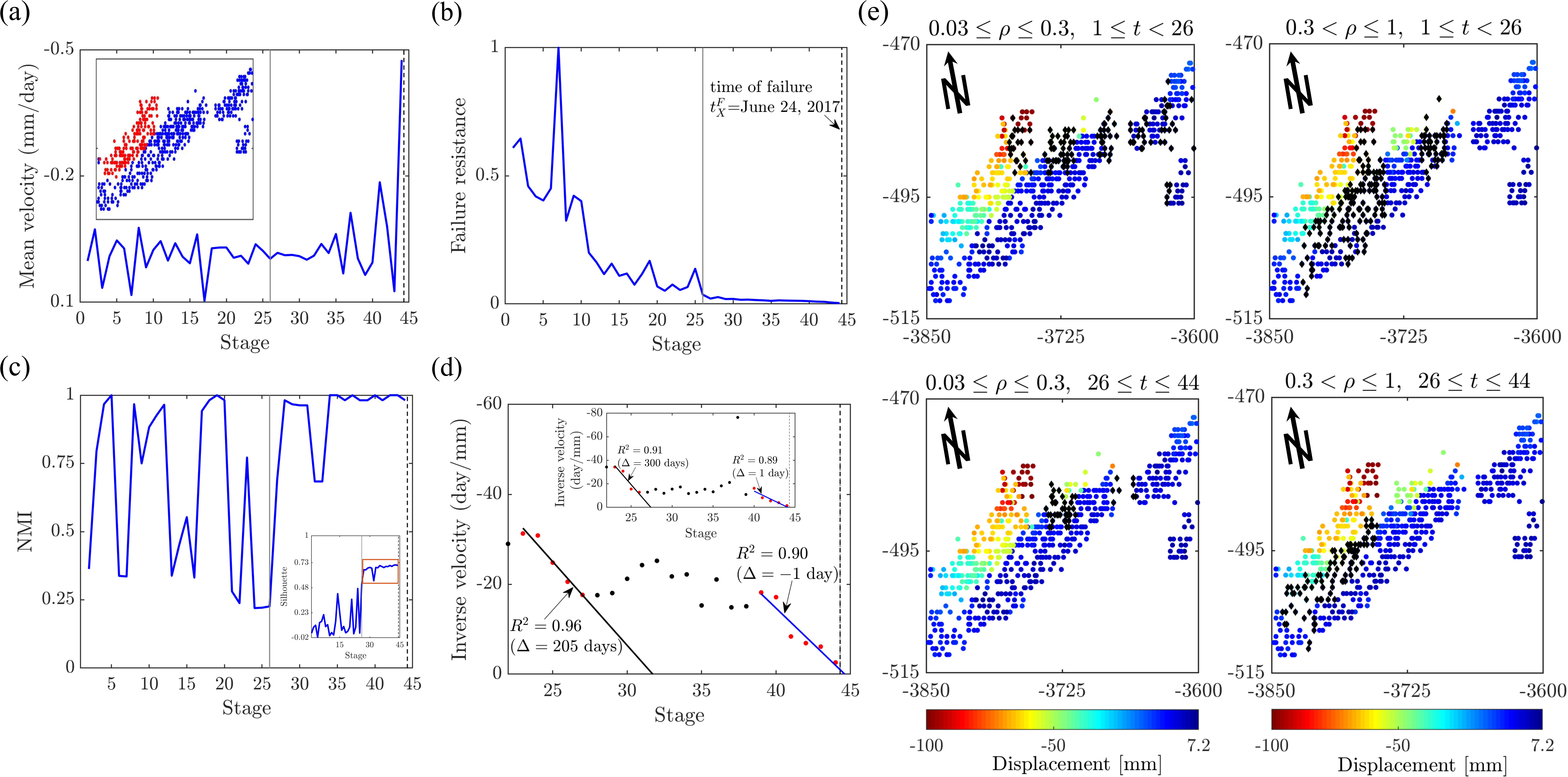}
\caption{(Color online) Slope Xinmo. Time evolution of: (a) the mean velocity of Xinmo with the failure location in red (inset); (b) the failure resistance $F^*$ from Step 2; (c) the NMI index and the global average silhouette score (inset); (d) the inverse mean velocity of $\Omega$ and of pixel $p$ (inset), over the time interval indicated by the red window in the inset of (c).  Vertical lines mark the regime change point $t^*_X$=26=August 23, 2016 (solid grey line) and time of failure $t^F_X=$June 24, 2017 (dashed black line).  (e) Cumulative predictions of the landslide boundary (black diamonds) are overlaid on top of the displacement map at $t^F_X$. Supplementary Movie Xinmo shows the evolution of $\Omega$ over the period of the monitoring campaign.}
\label{fig:xinmo}
\end{figure}

In all of the systems studied, SSSAFE uncovers three dynamical regimes over the course of the monitoring campaign, consistent with a compromise-in-competition between force and damage (Figures~\ref{fig:oppose} --\ref{fig:xinmo}).  In Biax, the global mean velocity steadily rises in PFR, before a sudden burst to a peak in the failure regime (Figure~\ref{fig:oppose} (a)).  Simultaneously, the opposite trend can be observed in the time evolution of the system's resistance to failure $F^*$, which decreases progressively as damage spreads in PFR, eventually dropping to its minimum value close to zero at stage 80 (Figure~\ref{fig:oppose} (b)).  Extensive published studies of this sample has shown that columnar force chains at stage 80 have lost considerable lateral support in the region of impending shear band, due to dilatancy~\cite{Tordesillas2007,TORDESILLAS2009706,TORDESILLAS2011265}.  While force redistributions around force chains continually occur during this period (recall Figure~\ref{fig:rewire}), ultimately, the degradation in the region precipitates collective force chain buckling at the peak stress ($t=98$, Figure~\ref{fig:oppose} (a) inset), culminating in a fully developed shear band at $t^F=104$, when kinematic clusters move in almost relative rigid-body motion (recall Figure~\ref{fig:pattern} (a)). These events were previously observed in various types of sand and photoelastic disk assemblies~\cite{tordesillas2015,PRE,TordesillasWalkerAndoViggiani2013}.  A consistent dynamics emerges in the time evolution of NMI and $S$ in Figure~\ref{fig:oppose} (c). From stage $t=80$, $S$ rises from around 0.5 before levelling off at the start of the failure regime at $t=104$; NMI stays close to 1 from $t=80$.  These trends imply a recurring bottleneck $B(\mathcal{F})$ as evident in $80 \le t \le 110$ of Figure~\ref{fig:oppose} (d), such that grains on either side progressively move collectively as one in opposite directions (Figure~\ref{fig:pattern} (a)).

\vspace{10mm}

At the field scale using radar data, SSSAFE delivers qualitatively similar trends for Mines 1 and 2 and Xinmo.  The presence of large fluctuations in Mines 1 and 2 (Figures~\ref{fig:mine1}-\ref{fig:mine2} (a-c)) is not surprising given these mines were operational with blasting, pumping, transport and drilling works taking place at various times over the course of the monitoring period.  Like in Biax, the failure resistance of Mine 1 drops close to zero well before failure (Figures~\ref{fig:mine1} (b)), with corresponding rises in NMI and S towards 1 (Figures~\ref{fig:mine1}(c)), as early as around $t^*_{1a}=69$, even though the rock slope appears intact with near-zero global mean velocity (Figures~\ref{fig:mine1} (a)).  This suggests that internal cracks and shear bands have started to propagate internally, as delineated by the preferred paths of damage (black points in Figures~\ref{fig:mine1} (d)).  Damage spread is to the extent that the capacity for force transfer between adjacent connected material points (pixels) along the {\it recurring} bottleneck $B(\mathcal{F})$ to the west, has significantly reduced even though there are still many remaining connections in the rock slope that keep it manifestly intact: Supplementary Movie Mine1 clearly shows that the west wall persists as an active area $\Omega$ from the beginning of the monitoring campaign.

As time advances towards failure, a second regime change point manifests: $t^*_{1b}=3322$=12:14 June 14 (Figures~\ref{fig:mine1} (b-c)).  During $69 \le t < 3322$, the interaction between the two regions of instability leads to an initial decline in $S$ while NMI stays close to 1 due to the persistence of the west wall cluster, the site that eventually collapses.   But the day before the collapse, $S$ sharply rises from $t^*_{1b}$.  This rise in $S$ suggests that the clustering pattern has now become incised in the slope to the extent that the clusters are now essentially undergoing relative motion along their common boundary, as $\Omega$ accelerates~\cite{TordesillasZhouBatterham2018,SDAT,KSAT,ZhouAOAD}.  This is corroborated by the INV analysis of $\Omega$ and of the fastest moving pixel $p$ in $\Omega$ (Figures~\ref{fig:mine1} (e)). The change point $t^*_{1b}$ improves on earlier work using a pattern mining approach which detects the time of imminent failure to be one to two hours later: $t$=14:53 June 14 \cite{SDAT} and $t$=13:16 June 14 \cite{KSAT}.

In Mine 1, multiple sites of instability interact mechanically in PFR.  Our method can reliably identify and differentiate these regions (Figure~\ref{fig:mine1}).  The west wall where catastrophic failure occurs can be distinguished early in PFR by the temporal persistence of the predicted landslide boundary (black points) in this area, in contrast to the eastern corner where this boundary only occasionally appears (Supplementary Movie Mine1). This intermittent dynamics in the latter is due to redundant force pathways which the system exploits to relieve the build up of stress in $B(\mathcal{F})$ along the west wall by diverting the forces and damage there to alternative paths, including to the competing slide to the east.   In laboratory samples undergoing quasi-brittle failure~\cite{ATSKCRMNJT2}, the interaction between competing cracks manifest in the form of stress redistributions along the preferred force pathways {\it between} the bottleneck where the macrocrack ultimately forms and the competing crack which later undergoes structural arrest (self-stabilize).  The same is observed in Mine 1: note the concentration of black points in the area between the actual failure region to the west and the competing slide to the east in $0.3 < \rho \le 1$, $69 \le t < 3322$ of Figure~\ref{fig:mine1} (d). This compromise-in-competition continues until all such paths are exhausted, $t = t^*_{1b} =3322$, from which time $B(\mathcal{F})$ remains fixed and becomes primed for uncontrolled crack propagation, along the landslide boundary ($3322 \le t \le 3600$, Figure~\ref{fig:mine1} (d)).  Mine 1 provides a good example of why early prediction of failure rests crucially on methods that can account for the spatiotemporal compromise-in-competition between force and damage pathways.  Essentially, the ultimate effect of stress redistributions is to delay failure, since any damage to $B(\mathcal{F})$ leads to a reduction in $F^*$. But there is an undesired concomitant which is the considerable uncertainty they pose for early prediction of failure, given damage is rerouted and concentrated elsewhere -- away from the region of impending failure in PFR~\cite{ATSKCRMNJT2}.

In Mine 2, a consolidated rock slope dominated by intact igneous rock that embodies many natural joints or faults.   There is the A-dominated stable regime over the first day of the monitoring period $1 \le t < 221$, where the global mean velocity fluctuates initially around 0, as $F^*$ portrays a decreasing trend (Figure~\ref{fig:mine2} (a-b)).  Trends in both NMI and Silhouette coefficient $S$ suggest that no substantial clustering structure in the kinematics developed on the first day: NMI fluctuates between 0 and 1 while S remains below 0.5 (Figure~\ref{fig:mine2} (c)). The system embodies redundant pathways to divert stresses away from the area of impending failure ($1 \leq t < 221$, Figure~\ref{fig:mine2} (d)).  In the A-B dominated mesoregime of PFR, $S$ progressively increases to 1, implying the emergence of collective motion ($221 \leq t \leq 1370$, Figure~\ref{fig:mine2} (d)).  As failure draws near, intracluster motions become coherent and near rigid-body, while intercluster motions become separated (Figure~\ref{fig:mine2} (c) inset), as the cluster corresponding to the location of impending failure $\Omega$ accelerates (Figure~\ref{fig:mine2} (e)).  We see these trends are precisely mirrored by the Normalized Mutual Information (NMI) of the clusters (Figure~\ref{fig:mine2} (c)). Note that the landslide boundary, shown at $t= 221$ in Figure~\ref{fig:mesofield} actually appears as early as $t= 104$ and persists up until $t=178$, which explains the high NMI scores. However, the kinematic clusters undergo a short period of change during $179 \leq t < 221$ which may reflect any number of perturbations on the mine site, including blasting.  Around the same time interval, large fluctuations can also be observed in $S$.  Close to and during the B-dominated stable regime, $S$ flattens out close to 1, indicative of a strong clustered motion.  Altogether, the evidence from $F^*$, $S$ and NMI marks a regime change point at $t^*_2$=221=13:39 August 20, which is just over 4 days prior to the collapse on $t^F_2$=1315=03:00 August 25. The INV analysis of $\Omega$ and the fastest moving pixel $p$ supports the progressive evolution to collapse at $t^F_2$(Figure~\ref{fig:mine2} (e)).  

For the Xinmo landslide, two key trends are evident from SSSAFE: (a) the regime change point on $t^*_{X} = 26$ = August 23, 2016 which is 10 months in advance of the actual time of collapse from when the active area $\Omega$ became fixed in the location that later became the rock avalanche source~\cite{Intrieri2017} (Figure~\ref{fig:xinmo} (a-e)); and (b) the development of cracks in a smaller competing failure zone above and to the east of the actual rock avalanche source (black diamonds in $0.03 \le \rho \le 0.3$, $1 \le t < 26$ of Figure~\ref{fig:xinmo} (e), Supplementary Movie Xinmo) as early as 2015.
Thus results here corroborate prior findings~\cite{Intrieri2017,Hu2018,fan2017,carla2019,WLF5} that the rock avalanche was generated from the upper part of the slope (Figure~\ref{fig:xinmo} (a) inset), near the mountain crest, and only afterwards did it entrain old landslide deposits along its way, which showed no sign of movement and played a role only in the propagation and not in the triggering phase.  In the months immediately following $t^*_X$, the velocities recorded were slow~\cite{CrudenVarnes1996} and predisposed to self-stabilization~\cite{hungr2014}.  

Interestingly, the shift in $\Omega$ from the eastern to the western flank of the source area match the reconstruction proposed by Hu et al.~\cite{Hu2018}, who attribute the triggering of the rock avalanche to an initial rockfall on June 24, 2017 at 05:39:07 (local time) and impacting the bedrock in the source area, where a crack network was present but still locked by persistent rock bridges.  While the long creep behavior of the source area demonstrates that the failure was the result of a long-term process, instead of the sudden outcome of an external impulse, it is credible that the trigger rockfall hypothesized by Hu et al.~\cite{Hu2018} could have imposed the final stresses needed to overcome the failure resistance (capacity $c$ in Equation \ref{eq:cap}) of the remaining connections in the recurring bottleneck $B(\mathcal{F}(t)), \, \forall t \ge t^*_{X}$.  That $B(\mathcal{F})$ persisted in the same location from August 23, 2016 strongly suggests a progressive degradation in rock strength all along this path, with the antecedent prolonged rainfall~\cite{fan2017} likely aiding this condition and rendering $B(\mathcal{F})$ increasingly poised for uncontrolled crack propagation in the lead up to the failure event on June 24, 2017 (Figure~\ref{fig:xinmo} (b)).  Finally, around 40 days prior to the collapse, SSSAFE’s predicted area of collapse began to manifest a linear trend in the temporal evolution of its inverse mean velocity, which delivered a time of failure $t^F_X$ a day later than the actual collapse.





SSSAFE offers a lead time of a day to weeks.  This constitutes sufficient forewarning to undertake evacuation and other response actions~\cite{guzzetti2020}.  SSSAFE takes only a few tens of seconds per time state to generate predictions on the likely region of failure: 30 seconds for Mine 1 and Xinmo, and 50 seconds for Mine 2, on a standard laptop computer with 8 cores 1.30 GHz CPU. Thus a prediction can be returned before the next measurement even for the most advanced radar technology (e.g., 1-5 minutes).  At this rate, a reasonable number of time states (e.g., 30 consecutive time states would take at most 30 minutes) to establish robustly the dynamics of the region of interest for the purposes of identifying $t^*$ and $t^{F}$.  

\section*{Conclusion}
\label{sec:conclude}

A holistic framework for Spatiotemporal Slope Stability Analytics for Failure Estimation (SSSAFE) is developed.  We demonstrate how SSSAFE can be applied to identify emergent kinematic clusters in the early stages of the precursory failure regime for four case studies of catastrophic failure: one at the laboratory scale using individual grain displacement data; and three slopes at the field scale, using line-of-sight displacement of a slope surface, from ground-based and space-borne radars.  The spatiotemporal dynamics of the kinematic clusters reliably predicts where and when catastrophic failure occurs. The clusters share a common boundary along the path of least failure resistance.  Here we found this path to precisely locate the impending, shear band in the laboratory sample and the landslide boundary in the natural and man-made slopes.  The regime change point is marked by intracluster (intercluster) motions becoming very similar or rigid-body (separated) which, in turn, induces a spatial pattern of physical partitions that become invariant in time through to failure.  Our findings illuminate a way forward to rationalize and refine decision-making from broad-area coverage monitoring data for improved geotechnical risk assessment and hazard mitigation.  To that end, ongoing efforts are focused on the extension of SSSAFE to a probabilistic platform~\cite{WangSS} that incorporates uncertainty systematically for various slopes and relevant scenario projections.


\section*{Appendix}
\label{sec:appendix}

{\bf A. The biaxial compression test data} 
The laboratory scale data set Biax is a well studied data set from a simulation of polydisperse spherical grains confined to planar biaxial compression test (Figure~\ref{fig:system})~\cite{Tordesillas2007}.  
It is governed by a classical DEM model~\cite{CundallStrack1979}, modified to incorporate a moment transfer to account for rolling resistance and thus capture the effects of non-idealized particle shapes.
A combination of Hooke's law, Coulomb's friction, and hysteresis damping is used to model the interactions between contacting particles.
An initially isotropic assembly of $5098$ spherical particles is prepared and subject to biaxial compression with motion constrained to the plane.
The ensemble is subjected to constant confining pressure, with a coefficient of rolling friction of $\mu^r = 0.02$.
The initial packing fraction is $0.858$.
A rolling resistance and a sliding resistance act at the contacts, both of which are governed by a spring up to a limiting Coulomb value of ${\mu}|{{\bf f}^n}|$ and ${\mu^r }{R_{min}}|{{\bf f}^n}|$, respectively, where ${{\bf f}^n}$ is the normal contact force and $R_{min}$ is the radius of the smaller of the two contacting particles.  A summary of all the interaction parameters governing the
contacts and other quantities relevant to this study is presented elsewhere \cite{Tordesillas2007}.

\bibliography{References}

\section*{Acknowledgments}
We thank our two anonymous reviewers whose insightful comments and suggestions helped improve and clarify this manuscript.  AT and RB thank Prof. Jinghai Li and Dr. Jianhua Chen for the stimulating discussions at the 2nd International Panel of Mesoscience. 
AT and SK acknowledge support from the U.S. Army International Technology Center Pacific (ITC-PAC) and US DoD High Performance Computing Modernization Program (HPCMP) under Contract No. FA5209-18-C-0002.

\section*{Author contributions statement}
AT conceived, designed and coordinated the research.  SK conducted the network flow experiments.  RB supplied the data for Mine 1.  LC and PB supplied the data for Mine 2. EI supplied the data for Xinmo.  AT and SK analyzed the results.  AT wrote the original draft and revisions.  All authors reviewed, and contributed to the writing of, the manuscript.

\section*{Additional information}

\noindent \textbf{Competing interests} \\
The authors declare no competing interests.

\noindent \textbf{Data availability}\\
The data that support the findings of this study are available from GroundProbe Pty Ltd but restrictions apply since these data were used under license for the current study.

\section*{Symbols and their meaning}

\begin{flushleft}
\begin{longtable}{p{.2\textwidth}  p{.7\textwidth}} 
\hline
\hline
\hline
Symbol & Meaning \\
\hline
$|\overrightarrow{\Delta u_{ij}}|$ & the magnitude of the relative displacement of two grains $i$ and $j$ (or two pixels) linked in $\mathcal{N}$\\
$\delta^+(v)$ & set of arcs leaving node $v$\\ 
$\delta^-(v)$ &  set of arcs entering node $v$\\ 
$\Gamma$ & a cut or virtual crack path in $\mathcal{N}$\\
$\Gamma_{min}$ & minimum cut in $\mathcal{N}$\\
$\mu^r$ & coefficient of rolling friction\\ 
$\rho(l)$ & cut ratio, the ratio of the number nodes in the smallest to largest components upon the removal of link $l$ in $\mathcal{T}$  \\
$\Omega$ & active and fastest-moving kinematic cluster, predicted location of failure\\

\hline
\hline
$A$ & set of arcs\\
$a(i)$ & average distance in the displacement state-space from $i$ to all other nodes in the same cluster\\
$B(\mathcal{F})$ & bottleneck, minimum cut of $\mathcal{F}$ \\
$b(i)$ &  average of the distances from $i$ to all the points in the other clusters\\
$c$ & capacity function \\
$c(e)$ & capacity of link $e$, threshold for damage, strength of contact, failure resistance of contact, proximity to failure \\
$c(\Gamma_{min})$ & capacity of the minimum cut $\Gamma_{min}$\\
$\vec{d}_\ell$ & displacement vector of grain (pixel) $\ell$\\
EWS & early warning systems\\
$e$ & link in $A$ corresponding to a contact in $\mathcal{N}$ \\
$\mathcal{F}$ & flow network\\ 
$F^*$ & failure resistance, global force flow capacity, bottleneck capacity\\
${{\bf f}^n}$ & normal contact force\\
$f(x)$ & total flow leaving the source node $q$\\
$G$ & directed network\\
$G^*$ & undirected link-capacitated network\\
$H(X(t))$ & entropy of the clustering assignment $X(t)$\\
$I(X(t);X(t-1))$ & mutual information between $X(t)$ and $X(t-1)$\\ 
$i, j$ & grain (pixel) in $\mathcal{N}$ which corresponds to a physical contact $e$\\
$k$ & sink node \\
$l$ & link in $\mathcal{T}$\\
$l_m$ & link of minimum weight on the path joining nodes $u$ and $v$ in $\mathcal{T}$\\
LOS & line-of-sight\\
$\ell$ & observation point (grain or pixel)\\
NMI & normalized mutual information\\
$\mathcal{N}$ & physical contact network\\
$n$ & number of nodes in $\mathcal{N}$, or $G^*$\\
$O$ & $(n-1) \times 3$ matrix whose entries are all zeros\\
$O[i,:]$ & $i^{th}$ row of $O$\\ 
$p$ & pixel with the highest rate of movement in $\Omega$ at the time of failure\\
PFR & precursory failure regime, mesoregime\\
$q$ & source node\\
$R_{min}$ & radius of the smaller of the two contacting particles\\
$S$ & Silhouette score for the whole body\\
$s(i)$ & Silhouette score for a node $i$\\
SSR & slope stability radar\\
$T$ & final time stage\\
$\mathcal{T}$ & Gomory-Hu tree\\
$t$ & time stage\\
$t^F, t_B^F, t_1^F, t_2^F, t_X^F$ & time of failure for a general case, Biax, Mines 1, 2 and Xinmo, respectively\\
$t^*, t_B^*, t_1^*, t_2^*, t_X^*$ & regime change point for a general case, Biax, Mines 1, 2 and Xinmo, respectively\\
$V$ & set of nodes\\
$v, u$ & node in $G^*$\\
$W,W'$ & two disjoint components of $\mathcal{N}$ corresponds to the removal of edges in $B(\mathcal{F})$ or removal of a link $l$ in $\mathcal{T}$\\
$w(l)$ & non-negative weight on link $l$ in $\mathcal{T}$\\
$w(l_m)$ & capacity of the minimum cut separating $u$ and $v$ in $\mathcal{N}$\\ 
$X(t)$ & clustering assignment at time stage $t$\\
$x(e)$ & flow of link $e$\\
 
\hline
\hline
\hline
\end{longtable}
\end{flushleft}

\end{document}